\documentclass[english]{article}
\usepackage[T1]{fontenc}
\usepackage[latin9]{inputenc}
\usepackage{amsmath}
\usepackage{amssymb}
\usepackage{stackrel}
\PassOptionsToPackage{normalem}{ulem}
\usepackage{ulem}

\makeatletter
\@ifundefined{date}{}{\date{}}
\makeatother

\usepackage{babel}
\begin{document}
\title{Numbers Extensions}
\author{David O. Zisselman}
\maketitle

\section{Introduction and abstract }

Over the course of the last 50 years, many questions in the field
of computability were left surprisingly unanswered. One example is
the question of $P$ vs $NP\cap co-NP$ . It could be phrased in loose
terms as ``If a person has the ability to verify a proof and a disproof
to a problem, does this person know a solution to that problem?''. 

When talking about people, one can of course see that the question
depends on the knowledge the specific person has on this problem.
Our main goal will be to extend this observation to formal models
of set theory $ZFC$: given a model $M$ and a specific problem $L$
in $NP\cap co-NP$ , we can show that the problem $L$ is in $P$
if we have ``knowledge'' of $L$.

In this paper, we'll define the concept of knowledge and elaborate
why it agrees with the intuitive concept of knowledge. Next we will
construct a model in which we have knowledge on many functions. From
the existence of that model, we will deduce that in any model with
a worldly cardinal we have knowledge on a broad class of functions. 

As a result we show that if we assume a worldly cardinal exists, then
the statement ``a given definable language which is provably in $NP\cap co-NP$
is also in $P$ `` is provable. 

Assuming a worldly cardinal , we show by a simple use of these theorems
that one can factor numbers in poly-logarithmic time. 

This article won't solve the $P$ vs $NP\cap co-NP$ question, but
its main result brings us one step closer to deciding that question. 

\section{Preliminaries }

Before I begin, and since the proofs use a few known theorems and
basic definitions, I shall quote the theorems which I'll use later.

\subsection{Definitions and notations \label{subsec:Definitions-and-notations}}
\begin{itemize}
\item Unless specifically mentioned otherwise $\mathbb{N}$ will denote
the set of natural numbers $\mathbb{N}=\left\{ 0,1,2,3,4,...\right\} $
which includes zero.
\begin{itemize}
\item As, strictly speaking, $\mathbb{N}$ can't be defined I will not use
the notation $\mathbb{N}$ except for intuitions. 
\end{itemize}
\end{itemize}
For the following definitions we will assume we have one universal
Turing machine (i.e a coding method of Turing machines) by which we
measure the length of other Turing machines 
\begin{itemize}
\item The Kolmogorov Complexity of a natural number $x$ is denoted by $\text{Kol}\left(x\right)$.
\end{itemize}

\paragraph{Definition (text to integer coding):}

Given a finite set of letters $\Sigma$ and string of text $T\in\Sigma^{*}$
one can code $T$ into a number in the following way:
\begin{enumerate}
\item fix a numbering on the set of letters
\begin{itemize}
\item for example if $\Sigma=\left\{ a,b,c,\ldots\right\} $ set $\#a=1,\#b=2,\ldots$ 
\end{itemize}
\item replace every letter in $T$ with its numbering, and get a number
of $\left|T\right|$ digits in base $\left|\Sigma\right|$
\begin{itemize}
\item for example if $\Sigma=\left\{ a,b,c,\ldots,z\right\} $ and $T="cbac"$
then the coding is
\begin{align*}
\#c+\#b\cdot26+\#a\cdot26^{2}+\#c\cdot26^{3} & =\\
2+1\cdot26+0\cdot26^{2}+2\cdot26^{3} & =17,604
\end{align*}
\end{itemize}
\end{enumerate}
Such a coding is called \textbf{text to integer coding.}

\subparagraph{Observation:}

text to integer coding can be preformed using \textbf{only} the arithmetic
operations addition, multiplication, exponentiation along with the
numbering on the set of letters (in the above example $\#a=0,\#b=1,\ldots$
)

\paragraph{Definition:}

The language of set theory is a single two place relation along with
the symbol of equality, that is $\left\{ \in,=\right\} $

\paragraph{Assumption:}

We will fix a numbering on the alphabet of first order logic along
with the language of set theory i.e a numbering for the alphabet 
\[
\left\{ "(",")","\neg","\forall","\exists","\wedge","\vee","\in","="\right\} 
\]

\paragraph{Notation:}

A \textbf{Turing machine} will refer to a RAM computational machine. 

\paragraph{Definition:}

A \textbf{structure }for the language of set theory is a set $A$
along with a two place relation $\left(\in\right)\subset A^{2}$ . 

I.e. a set of elements $A$ along with a subset $R\subset A^{2}$
which we will denote by $a\in b\Leftrightarrow aRb$ 

\paragraph{Definition ($ZF$):}

the set of axioms on the language of set theory called $ZF$ (Zermelo--Fraenkel
axiomatic set theory) is the following set:
\begin{enumerate}
\item Axiom of extensionality $\forall x\forall y\left(\left(\forall z\left(z\in x\leftrightarrow z\in y\right)\rightarrow x=y\right)\right)$
\item Axiom of regularity $\forall x\left(\exists a\left(a\in x\right)\rightarrow\exists y\left(y\in x\wedge\neg\exists z\left(z\in y\wedge z\in x\right)\right)\right)$
\item Axiom schema of specification $\forall z\forall w_{1}\forall w_{2}...\forall w_{n}\exists y\forall x\left(x\in y\leftrightarrow\left(x\in z\wedge\phi\left(w_{1},...,w_{n},x\right)\right)\right)$
\item Axiom of pairing $\forall x\forall y\exists z\left(x\in z\wedge y\in z\right)$
\item Axiom of union $\forall F\exists A\forall Y\forall x\left(\left(x\in Y\wedge Y\in F\right)\rightarrow x\in A\right)$
\item Axiom schema of replacement 
\[
\!\!\!\!\!\!\!\!\!\!\!\!\!\!\!\!\!\!\!\!\!\!\!\!\!\!\!\!\!\!\!\!\!\!\!\forall A\forall w_{1}\forall w_{2}...\forall w_{n}\left[\forall x\left(x\in A\rightarrow\exists!y\phi\left(w_{1},...,w_{n},x,y\right)\right)\rightarrow\exists B\forall x\left(x\in A\rightarrow\exists y\left(y\in B\wedge\phi\left(w_{1},...,w_{n},x,y\right)\right)\right)\right]
\]
\item Axiom of infinity $\exists X\left(\emptyset\in X\wedge\forall y\left(y\in X\rightarrow S\left(y\right)\in X\right)\right)$
\item Axiom of power set $\forall x\exists y\forall z\left(z\subseteq x\rightarrow x\in y\right)$
\end{enumerate}
\begin{itemize}
\item Definition: The axiom of choice (AC) is the following statement:
\[
\forall X\left[\emptyset\not\in X\Rightarrow\exists f:X\rightarrow\cup X\,\,\,\forall a\in X\left(f\left(a\right)\in a\right)\right]
\]
\item Definition: the set of axioms $ZF$ with AC is called $ZFC$ axioms
set.
\item Definition: the set of axioms $ZF$ without regularity is called $ZF^{-}$ 
\end{itemize}
Unless mentioned otherwise - a language will refer only to a definable
language.

\paragraph{Definition:}

A model of ZF/ZFC/ $ZF^{-}$ etc is a structure for the language of
set theory which satisfies the appropriate axioms of ZF/ZFC/ $ZF^{-}$
respectively. 

\paragraph{Definition and notation of $\omega$:}

Within a model $M$ of $ZF$ the standard way to construct the ``natural
numbers'' within the model is to define 0 to be $\emptyset$ and
$S\left(x\right)=x\cup\left\{ x\right\} $ (successor operation).
So $1=\left\{ \emptyset\right\} $ and $2=\left\{ \left\{ \emptyset\right\} ,\emptyset\right\} $
and so on. The minimal set which contains $\emptyset$ and is closed
under the operation $S$ is called $\omega$ and the its existence
is guaranteed by the axiom of infinity. However, and even though we
conceive $\omega$ as $\mathbb{N}$, they may be very different objects
in some models. If in a certain model , the set $\omega$ differs
from our regular notion of natural numbers, we can consider such a
model as having ``non-standard'' arithmetic.
\begin{itemize}
\item The notation $\omega$ will be a more suitable notation than the imprecise
notation $\mathbb{N}$.
\end{itemize}

\subparagraph{Notation:}
\begin{itemize}
\item Given a model $M$ of $ZF$ denote $\omega\left(M\right)$ to be the
$\omega$ of the model $M$.
\item Given an $n\in\omega\left(M\right)$ denote 
\[
\left[n\right]\triangleq\left\{ i\in\omega\left(M\right)|i\leq n\right\} =\left\{ 0,1,2,\ldots,n\right\} 
\]
\end{itemize}

\paragraph{Definitions of consistencies:}
\begin{itemize}
\item for a given $j\in\omega$, $j-con\left(ZF\right)$ is defined by the
following statement: \\
Given a model $M$ of $ZF$, one can find a sequence $M_{1},M_{2},M_{3},...,M_{j}$
within $M$ such that $M_{1}$ is a model of $ZF$ and a set in $M$
and for every $i<j$ $M_{i+1}$ is a set of $M_{i}$ and a model of
$ZF$ and $\in_{i+1}$ (as a subset of $M_{i+1}\times M_{i+1}$) is
a set of $M_{i}$.
\item the statement $\omega-con\left(ZF\right)$ is the following statement:
\\
Given a model $M$ of $ZF$ for all $j\in\omega\left(M\right)$, one
can find a sequence of models of $ZF$ $M_{1},M_{2},M_{3},...,M_{j}$
within $M$ such that $M_{1}$ is a set in $M$ and for every $i<j$
$M_{i+1}$ is a set of $M_{i}$ and $\in_{i+1}$ (as a subset of $M_{i+1}\times M_{i+1}$)
is a set of $M_{i}$.
\item the statement $\left(\omega+1\right)-con\left(ZF\right)$ is the following
statement: \\
Given a model $M$ of $ZF$, one can find a model $M_{0}$ of $ZF$
within $M$ s.t $M_{0}$ satisfies $\omega-con\left(ZF\right)$ .
\item The same definitions $j-con\left(ZFC\right)$ apply for $j\in\omega$
or $j=\omega$ or $j=\omega+1$.
\end{itemize}

\subparagraph{Observation:}
\begin{itemize}
\item For a model $M$ of $ZF$ and for $j\in\omega\left(M\right)$ one
can define the formula $j-con\left(ZFC\right)$ within the model $M$
arithmetic. This applies to the case where $M$ has ``non-standard''
arithmetic and $j$ is a non-standard number as well. The same holds
for $j-con\left(T\right)$ and for effective theories $T$.
\end{itemize}

\paragraph{Definition (Von Neumann universe): }

Within a model $V,\in$ of ZFC define:
\begin{itemize}
\item $V_{0}=\emptyset$
\item $V_{\alpha+1}=P\left(V_{\alpha}\right)$ where $P\left(X\right)$
is the power set of $X$
\item $V_{\alpha}=\cup_{\beta<\alpha}V_{\beta}$ for $a=\cup_{\beta<\alpha}\beta$
limit ordinal.
\end{itemize}

\paragraph{Definition (transitive set):}

Within a model $V,\in_{V}$ of ZF a set $A$ is called \textbf{transitive}
if for every $x,y$ sets of $V$ if $x\in_{V}A$ and $y\in_{V}x$
then $y\in_{V}A$.

\subparagraph{Definition (set model):}

Given a ZF model $V_{1},\in_{V_{1}}$ another model $V_{2},\in_{V_{2}}$
of ZF is said to be \textbf{inside} $V_{1}$ (or a set of $V_{1}$)
if the following hold:
\begin{itemize}
\item $V_{2}$ is a set within $V_{1}$ 
\item all sets of $V_{2}$ are also sets in $V_{1}$.
\item $\in_{V_{2}}$ as a set of 2-topuls (i.e as a subset of $\left(V_{2}\right)^{2}$)
is a set of $V_{1}$
\end{itemize}

\paragraph{Lemma:}

Let $V_{2}$ and $V_{1}$ be $ZF$ models if $V_{2}$ is a set of
$V_{1}$, $a$ is a set of $V_{2}$, $b\in_{V_{2}}a$, then $b$ is
a set of $V_{2}$ and $b$ is a set of $V_{1}$.

\paragraph{Proof:}

As $\in_{V_{2}}$ is a subset of $\left(V_{2}\right)^{2}$, $b\in_{V_{2}}a$
means $b$ is a set of $V_{2}$. As $V_{2}$ is inside $V_{1}$, this
means $b$ is also a set of $V_{1}$. $\blacksquare$

\subsubsection{Worldly cardinal}

\paragraph{Definition:}

A cardinal $k$ is called \textbf{worldly }if\textbf{ $\left(V_{k},\in\right)$
}is a $ZFC$ model.

\subsubsection{Mostowski\textquoteright s Collapsing Theorem \label{subsec:Mostowski=002019s-Collapsing-Theorem}}

\paragraph{Theorem:}
\begin{enumerate}
\item If $E$ is a well-founded and extensional relation on a class $P$,
then there is a transitive class $M$ and an isomorphism $\pi$ between
$\left(P,E\right)$ and $\left(M,\in\right)$. The transitive class
$M$ and the isomorphism $\pi$ are unique.
\item In particular, every extensional class $P$ is isomorphic to a transitive
class $M$. The transitive class $M$ and the isomorphism $\pi$ are
unique.
\item In case (2), if $T\subset P$ is transitive, then $\pi x=x$ for every
$x\in T$.
\end{enumerate}

\paragraph{Proof reference:}

Please refer to \cite{Kunen - set theory} Chapter 6 ``The Axiom
of Regularity'' (pg. 69) theorem 6.15.

\subsubsection{Lowenheim Skolem theorem \label{subsec:Lowenheim-Skolem-theorem}}

\paragraph{Theorem:}

Every infinite model for a countable language has a countable elementary
sub model.

\paragraph{Proof reference:}

Please refer to \cite{Kunen - set theory} Chapter 12 ``Models of
Set Theory'' (pg. 157) theorem 12.1.

\subsection{Forcing\label{subsec:Forcing}}

\paragraph{Definition:}

A partial order set (\textbf{POS}) is a triple $\left(P,\leq_{P},0_{P}\right)$
s.t $P$ is a set, $\leq_{P}$ is a partial order on $P$ and $0_{P}$
is a minimal element in $P$.

Let $\left(P,\leq_{P},0_{P}\right)$ be a POS.  Define:
\begin{itemize}
\item a set $D$ is \textbf{dense} with respect to $P$ if 
\[
\forall p\in P\exists q\in D\,\,q\leq_{P}p
\]
\item $G\subset P$ is called a \textbf{filter on $P$} if it satisfies
the following three conditions::
\begin{itemize}
\item $0_{P}\in G$
\item $\forall p,q\in P\,\,\left(q\leq_{P}p\,\,\wedge\,\,p\in G\right)\rightarrow q\in G$
\item $\forall p,q\in G\,\,\exists r\in G\,\,r\geq_{P}p\,\,\wedge\,\,r\geq_{P}q$
\end{itemize}
\item For a collection of sets $M$ (which may be a model of $ZF$) and
a filter $G$ on $P$, we say that $G$ \textbf{generic over }$M$
if for every dense set $D\in M$ we have $G\cap D\neq\emptyset$.
\item Given $E\subset P$ and $p\in P$, we say that $E$ is \textbf{dense
above $p$ }if 
\[
\forall q\geq p\,\,\exists r\geq q\,\,\left(r\in E\right)
\]
\item $\tau$ is a\textbf{ $P$ name} if $\tau$ is a relation and 
\[
\forall\left\langle \sigma,p\right\rangle \in\tau\,\left[\tau\text{ is a \ensuremath{P} name \ensuremath{\wedge\,p\in P}}\right]
\]
\item For $M$ a model of $ZF$ the $P$ \textbf{names} \textbf{in} $M$
are 
\[
M^{P}=\left\{ \tau\in M\,|\,\tau\text{ is a \ensuremath{P} name in \ensuremath{M}}\right\} 
\]
\item For $M$ a model of $ZF$ and $G$ a filter on $P$, the\textbf{ valuation
of a name }is
\[
Val\left(\tau,G\right)=\tau_{G}=\left\{ Val\left(\sigma,G\right)\,|\,\exists p\in G\,\,\left\langle \sigma,p\right\rangle \in\tau\right\} 
\]
\item For $M$ a model of $ZF$ and $G$ a filter on $P$ define 
\[
M\left[G\right]=\left\{ \tau_{G}\,|\,\tau\in M^{P}\right\} 
\]
\end{itemize}

\subparagraph{Remark:}
\begin{itemize}
\item $P$ names and valuations are both defined recursively. 
\end{itemize}

\subparagraph{Notation:}
\begin{itemize}
\item for $p,r\in P$ denote $p\bot r$ if $\neg\exists q\in P\,\,\left(q\leq p\,\,\wedge r\leq p\right)$
\end{itemize}

\paragraph{Theorems:}
\begin{enumerate}
\item Let $P$ be a POS and $M$ a countable collection of sets and let
$p\in P$. Than there exists a generic filter $G$ over $M$ s.t $p\in G$.
\item If $M$ is a transitive model of $ZFC$ and $G$ is a generic filter
over $M$ and $P\in M$ is a POS s.t
\[
\forall p\in P\,\,\exists q,r\in P\,\,\left(p\leq q\,\,\wedge\,\,p\leq r\,\,\wedge\,\,q\bot r\right)
\]
 then $G\not\in M$.
\item Given a countable transitive model $M$ of $ZFC$ and $P\in M$ and
a a generic filter $G$ over $M$ and $p\in G$ and $E\in M$ . If
$E\subset P$ is such that $E$ is dense above $p$, then $G\cap E\neq\emptyset$.
\item If $M$ is a transitive model of $ZFC$ and $P\in M$ and $G$ is
a generic filter over $M$, then $G\in M\left[G\right]$.
\item If $M$ is a countable transitive model of $ZFC$ and $P\in M$ and
$G$ is a generic filter over $M$, then $M\left[G\right]$ is a countable
transitive model of $ZFC$.
\end{enumerate}

\paragraph{Proof references:}

Please refer to \cite{Kunen - set theory} Chapter 7 ``Forcing'':
\begin{enumerate}
\item pg. 186 Lemma 2.3.
\item pg. 187 Lemma 2.4.
\item pg. 192 Lemma 2.19 $\left(ii\right)$.
\item pg. 190 Lemma 2.13
\item pg. 201 Theorem 4.2 for the fact that $M\left[G\right]$ holds $ZFC$.
$M\left[G\right]$ is countable (as a countable union of countable
sets) and transitive by it's construction.
\end{enumerate}

\subsection{Basic Theorems}

\subsubsection{Chaitin's incompleteness theorem (1971) \label{subsec:Chaitin's-incompleteness-theorem} }

\paragraph{Theorem:}

Let $V$ be a model of $ZF$. Let $T$ be an effective consistent
set of axioms. Then, there exist $L\in\mathbb{\omega}\left(V\right)$
(which depends on the set of axioms) such that for every $x\in\mathbb{\omega}\left(V\right)$
the statement $L\leq Kol\left(x\right)$ can't be proven from the
set $T$ .

For completeness of the article I'll add proof to this theorem:

The reader may want to consider the case of $\omega\left(V\right)=\mathbb{N}$
at first.

\paragraph{Proof: }

Within $V$, any proof of claim $\phi$ from $T$ is a number within
$\mathbb{\omega}\left(V\right)$. Denote $w\in\omega\left(V\right)$
to be the first proof of a claim in $\left\{ L\leq\text{Kol}\left(x\right)|x\in\omega\left(V\right)\right\} $
and let $x'\in\mathbb{\omega}\left(V\right)$ be the number s.t $T\underset{w}{\vdash}L\leq\text{Kol}\left(x'\right)$
(the proof $w$ proves that $L\leq\text{Kol}\left(x'\right)$). Create
a Turing machine $M$ that goes over $y\in\mathbb{\omega}\left(V\right)$
(in the regular order) and checks if $y$ is a proof of the statement
$L\leq\text{Kol}\left(x\right)$ and halts if it is and prints $x$.
The size of the TM $M$ (as TM could be coded as natural numbers too)
is $\log\left(L\right)+C$ for a fixed $C\in\mathbb{\omega}\left(V\right)$.
As, the decimal representation of $L$ is $\log\left(L\right)$ and
$C$ is the extra size required to represent the theory $T$ and the
operation of $M$. 

The TM $M$ is a TM that prints $x'$ and its size is $\log\left(L\right)+C$
and as $L\leq\text{Kol}\left(x'\right)$ it follows that $L\leq\log\left(L\right)+C$
(as we've just shown a TM of size $\log\left(L\right)+C$). the inequality
$L\leq\log\left(L\right)+C$ can't hold for a sufficiently large $L$. 

\subparagraph{Notation:}

given an effective consistent set of axioms $T$ define $\mathcal{L}\left(T\right)$
to be the minimal number s.t $L>\log\left(L\right)+C$. 

\paragraph{Remarks:}
\begin{itemize}
\item Notice that $\mathcal{L}\left(T\right)\leq2C$.
\item The notation $\mathcal{L}\left(T\right)$ doesn't mention $V$ because
the value doesn't depend on $V$ to a large extent (and it will later
be stated and proved).
\item here we abuse the notation somewhat, as $T$ is not a set of axioms,
but also a coding of a TM that identifies said set axioms and as such
the value of $\mathcal{L}\left(T\right)$ as defined above depends
on the representation of the machine.
\end{itemize}

\subsubsection{Compactness theorem \label{subsec:Compactness-theorem}}

\subparagraph{Theorem:}

Assume $T$ is a set of axioms such that every finite subset $A\subset T$
, $A$ is consistent, then $T$ is consistent.

\subparagraph{Proof reference:}

Please refer to \cite{Notes on logic and set theory - Johnstone}

Corollary 3.8 in chapter 3 pg. 27.

\subsubsection{Binary tree construction \label{subsec:binary-tree-construction}}

\paragraph{Theorem:}

Let $M$ be a model of $ZF$ let $\check{z}\in\omega\left(M\right)$
be a fixed number, $a\in\omega\left(M\right)^{\left[\check{z}\right]}$
is any sequence of numbers of length $\check{z}$ within the model. 

Then exist $t\in\omega\left(M\right)$ , which encodes a Turing machine
that given $k\leq\check{z}$ in binary representation calculates $a_{k}$
in $\left\lceil \log_{2}\left(k+1\right)+\log_{2}\left(a_{k}\right)+1\right\rceil $
computational steps.

\paragraph{Intuition:}

The key idea here is to build a binary tree with the values of $a_{k}$
and create a Turing machine which fetch $a_{k}$ using that binary
tree. Then the access time is the size of the representation of the
number $k$ which is $\left\lceil \log_{2}\left(k+1\right)+1\right\rceil $
and along with the time to write the result is $\left\lceil \log_{2}\left(k+1\right)+\log_{2}\left(a_{k}\right)+1\right\rceil $.

We will assume w.o.l.g that traveling through one of the edges of
the tree takes one computation step

\paragraph{Proof:}

Let $z\in\omega\left(M\right)$ be such that $\check{z}\leq2^{z}$.
Create the following data structure:
\begin{itemize}
\item a binary tree where every node has a value (in $\omega\left(M\right)$)
\item The depth of the tree is $z$.
\begin{itemize}
\item the value of $a_{0}$ will be held in a special place in memory.
\end{itemize}
\item The root will hold the value of $a_{1}$ (a number in $\omega\left(M\right)$)
\item by induction, if a node $k$ held the value of $a_{n}$ then:
\begin{itemize}
\item The left son of $k$ will hold the value of $a_{2n}$
\item The right son of $k$ will hold the value of $a_{2n+1}$
\end{itemize}
\end{itemize}
The Turing machine will go over the data structure from the root based
on the binary representation of $k$ (right to left), for every $1$
it would go to the right son and for $0$ it would go to the left
son.

For example, if we want to fetch $a_{13}$ as $13=1101_{2}$ then
for the root we would go right then right then left and finally right
to get to $a_{13}$. We will assume w.o.l.g that traveling through
one of the edges of the tree takes one computation step (because,
for example, such an operation is implemented in the hardware) and
thus this process takes the size of the representation of $k$ and
the time to write the data $a_{k}$ which is at most $\left\lceil \log_{2}\left(k+1\right)+\log_{2}\left(a_{k}\right)+1\right\rceil $
steps. $\blacksquare$

Please note that using this construction as $\check{z}$ gets larger
and larger so does $t$. 

\section{Extension and Lattices - first part.\label{sec:Lattices}}

On this part I'll describe the principle of numbers extensions and
show that given a definable function ,which satisfies certain conditions,
one can build a sequence of arithmetic models in which a sequence
of bounded Turing machines exist that calculates this function in
logarithmic time to an increasingly larger numbers (i.e a model of
$ZFC$ exist such that it's $\omega$ has a Turing machine that calculates
in logarithmic time...)

\subsubsection{Definable function\label{subsec:Definable-function}}

\paragraph{Definition: \label{par:Definition definable language}}

Assume $\phi\left(n_{1},n_{2}\right)$ is a two variable formula in
$ZF$ language and let $V$ be a model of $ZF$. We say that a set
$A^{2}\subset\omega\left(V\right)^{2}$ in model $V$ is definable
using formula $\phi$ if these conditions hold 
\[
A=\left\{ \left(x,y\right)\in\omega\left(V\right)^{2}|\phi\left(x,y\right)\right\} 
\]

\paragraph{Examples:}
\begin{enumerate}
\item The set of numbers $n,k$ s.t $n$ is divisible by $k$ is definable
using the formula $\phi\left(a,b\right)=\left(\exists c\in\omega\right)\left(b\cdot c=a\right)$.
\item Given a lexicographic coding of 3-sat formulas. Denote the formula
$s_{3}sat\left(n,k\right)$ ``the formula $\#n$ is a 3-sat formula
that can be satisfy using the assignment $\#k$''
\item Given a coding of Turing machines, the set of Turing machines which
halts by the $k$ step is definable using the formula $\phi\left(n,k\right)$
``$n$ is a TM that halts by the $k$ step''.
\end{enumerate}

\paragraph{Definition:}

Within a ZF model $V$. A two variable formula $\phi\left(n,k\right)$
defines a \textbf{function }in\textbf{ $V$} if 
\[
V\vDash\left(\forall x\in\omega\right)\left(\exists!y\in\omega\right)\phi\left(x,y\right)
\]

\paragraph{Definition:}

A two variable formula $\phi\left(n,k\right)$ defines a \textbf{function}
in $ZF$ if 
\[
ZF\vdash\left(\forall x\in\omega\right)\left(\exists!y\in\omega\right)\phi\left(x,y\right)
\]

\subsubsection{Definition of model extension \label{subsec:model-extension-definition}}

\paragraph{Definition: }

Given two models of $ZF$ , $V_{1},V_{2}$ we say that $V_{2}$ number
extends $V_{1}$ using mapping $f:\omega\left(V_{1}\right)\rightarrow\omega\left(V_{2}\right)$
if the arithmetic operations of $V_{1},V_{2}$ i,e 
\[
+_{1},+_{2},\cdot_{1},\cdot_{2}\wedge_{1},\wedge_{2}
\]
accordingly are respected by $f$. Meaning: 
\begin{enumerate}
\item $\forall a,b\in\omega\left(V_{1}\right)\left(f\left(a+_{1}b\right)=f\left(a\right)+_{2}f\left(b\right)\right)$
\item $\forall a,b\in\omega\left(V_{1}\right)\left(f\left(a\cdot_{1}b\right)=f\left(a\right)\cdot_{2}f\left(b\right)\right)$
\item $\forall a,b\in\omega\left(V_{1}\right)\left(f\left(a^{b}\right)=f\left(a\right)^{f\left(b\right)}\right)$
\item $f\left(0_{V_{1}}\right)=0_{V_{2}}$and $f\left(1_{V_{1}}\right)=1_{V_{2}}$
\end{enumerate}
I'll denote it by $V_{1}\overset{f}{\longrightarrow}V_{2}$. 

\subparagraph{Remark:}
\begin{itemize}
\item Please note that there is no assumption that inductive arguments on
elements of $\omega\left(V_{1}\right)$ are ``transferred'' (in
some way or another) into $\omega\left(V_{2}\right)\cap f\left(\omega\left(V_{1}\right)\right)$.
\item The function $f$ itself isn't assumed to be a set of either $V_{1}$
or $V_{2}$.
\end{itemize}

\subsubsection{Definition of model initial segment extension}

\paragraph{Definition:}

If $V_{1}\overset{f}{\longrightarrow}V_{2}$. and the extension satisfies
\[
\forall b\in\omega\left(V_{1}\right)\forall a\in\omega\left(V_{2}\right)\exists c\in\omega\left(V_{1}\right)\left(a<f\left(b\right)\rightarrow a=f\left(c\right)\right)
\]
 the extension is said to be \textbf{initial segment extension} (i.s
extension) and will be denoted by $V_{1}\stackrel[I]{f}{\rightarrow}V_{2}$
.

\paragraph*{Please note:}
\begin{enumerate}
\item The intuitive meaning of i.s. extension is that the model $V_{2}$
has ``more numbers'' than $V_{1}$ but the numbers $V_{1}$ ``behave
the same'' on $V_{2}$ ``lower parts''.
\item The function $f$ is usually undefinable from the models themselves.
\item The fact that $V_{1}\overset{f}{\rightarrow}V_{2}$ for two models
doesn't imply that the two models are elementary equivalent. Nor does
it mean that they agree on $\omega$ attributed formulas.
\end{enumerate}

\subsubsection{Definition of the propriety condition \label{propriety-condition-definition}}

\paragraph{Definition:}

Given a two variable formula $\phi\left(n,k\right)$ (in the language
if set theory) we say that the set defined by $\phi$ satisfies the
\textbf{propriety} condition if for every two models $V_{1},V_{2}$
of $ZF$ such that $V_{2}$ is a set of $V_{1}$ (and thus $V_{1}\underset{I}{\overset{f}{\longrightarrow}}V_{2}$)
it holds true that 
\[
\left(\forall n,k\in\omega\left(V_{1}\right)\right)\left(V_{1}\vDash\phi\left(n,k\right)\Leftrightarrow V_{2}\vDash\phi\left(f\left(n\right),f\left(k\right)\right)\right)
\]
.

\paragraph{Please note:}
\begin{itemize}
\item The meaning of the propriety condition is that one can determine if
a certain $a,b$ is a member of the set only by looking at arithmetic
operation on it an initial segment of the model. 
\begin{itemize}
\item One may think of the propriety condition of a function as a ``uniformly
recursive'' function 
\end{itemize}
\end{itemize}

\subsection{Basic properties of extension\label{subsec:basic-properties-of-extension}}

\paragraph{Lemma:}

Given two models of $ZF$ , $V_{1},V_{2}$ s.t 
\begin{itemize}
\item $V_{2}$ is a set within $V_{1}$ 
\item $\in_{2}$ is a set (of pairs) within $V_{1}$
\end{itemize}
then a mapping $f\in V_{1}$, $f:\omega\left(V_{1}\right)\rightarrow\omega\left(V_{2}\right)$
exist s.t $V_{1}\underset{I}{\overset{f}{\longrightarrow}}V_{2}$.

\paragraph{Proof:}

The construction of $f$ is by induction:
\begin{itemize}
\item The base case: define $f\left(0_{v_{1}}\right)=0_{V_{2}}$.
\item For $a=b+1$ where $a,b\in_{V_{1}}\omega\left(V_{1}\right)$ define
$f\left(a\right)=f\left(b+_{1}1\right)=f\left(b\right)+_{V_{2}}1$
\item As every natural number is either $0$ or $b+1$ for another natural
number $b$, $f$ is defined.
\end{itemize}
Recall that the definitions of $+,\times,\wedge$ are inductive as
well. For $a,c\in_{V_{1}}\omega\left(V_{1}\right)$ 
\begin{itemize}
\item If $c=0_{V_{1}}$ then by definition:
\begin{itemize}
\item $a+c=a+0=a$
\item $a\cdot c=a\cdot0=0$
\item $a^{c}=a^{0}=1$ and if $a\neq0$ then $c^{a}=0^{a}=0$
\end{itemize}
\item Therefore if $c=0$, $f$ holds the equalities in \ref{subsec:model-extension-definition}:
\begin{itemize}
\item $f\left(a+_{1}c\right)=f\left(a+_{1}0\right)=f\left(a\right)=f\left(a\right)+_{2}0_{V_{2}}$.
\item $f\left(a\cdot_{1}c\right)=f\left(a\cdot_{1}0_{V_{!}}\right)=f\left(0_{V_{1}}\right)=0_{V_{2}}=f\left(a\right)\cdot_{2}0_{V_{2}}=f\left(a\right)\cdot_{2}f\left(c\right)$
\item $f\left(a^{c}\right)=f\left(a^{0}\right)=f\left(1_{V_{1}}\right)=1_{V_{2}}=f\left(a\right)^{0_{V_{2}}}=f\left(a\right)^{f\left(c\right)}$
\item If $a\neq0$ then $f\left(c^{a}\right)=f\left(0^{a}\right)=f\left(0_{V_{1}}\right)=0_{V_{2}}=f\left(c\right)^{f\left(a\right)}$
\end{itemize}
\item For $c\neq0$ exist $b\in_{V_{1}}\omega\left(V_{1}\right)$ s.t $c=b+1$
and thus by definition:
\begin{itemize}
\item $a+c=a+\left(b+1\right)=\left(a+1\right)+b$
\item $a\cdot c=a\cdot\left(b+1\right)=\left(a\cdot b\right)+a$
\item $a^{c}=a^{b+1}=a\cdot\left(a^{b}\right)$
\end{itemize}
\item Therefore if $c\neq0$, $f$ holds the equalities in \ref{subsec:model-extension-definition}:
\begin{itemize}
\item $f\left(a+_{1}c\right)=f\left(\left(a+_{1}1\right)+_{1}b\right)=\left(f\left(a\right)+_{2}1_{V_{2}}\right)+_{2}f\left(b\right)$.
\item $f\left(a\cdot_{1}c\right)=f\left(\left(a\cdot_{1}b\right)+_{1}a\right)=f\left(a\cdot b\right)+_{2}f\left(a\right)=f\left(a\right)\cdot_{2}f\left(b\right)+_{2}f\left(a\right)=f\left(a\right)\left(f\left(b\right)+_{2}1_{V_{2}}\right)=f\left(a\right)f\left(b+_{1}1\right)=f\left(a\right)\cdot_{2}f\left(c\right)$
\item $f\left(a^{c}\right)=f\left(a\cdot_{1}\left(a^{b}\right)\right)=f\left(a\right)\cdot_{2}f\left(a^{b}\right)=f\left(a\right)\cdot_{2}f\left(a\right)^{f\left(b\right)}=f\left(a\right)^{f\left(b\right)+1_{V_{2}}}=f\left(a\right)^{f\left(c\right)}$
. 
\end{itemize}
\end{itemize}
So far we've seen that $V_{1}\overset{f}{\rightarrow}V_{2}$. 

This $f$ also holds the extension property as well: 

As $V_{2}$ is a set of $V_{1}$, in $V_{1}$ one can define the following
set 
\[
S=\left\{ b\in\omega\left(V_{1}\right)\,|\,\exists a\in\omega\left(V_{2}\right)\forall c\in\omega\left(V_{1}\right)\left(\left(a<_{2}f\left(b\right)\right)\wedge\left(a\neq f\left(c\right)\right)\right)\right\} 
\]
if $S$ isn't empty and as $S$ is a set of natural numbers in $\omega\left(V_{1}\right)$
it has a minimum. Denote 
\[
s=\min S
\]
 notice that $s\neq0$ as $f\left(0_{V_{1}}\right)=0_{V_{2}}$ and
$\neg\exists a\in\omega\left(V_{2}\right)\,\,\left(a<0_{V_{2}}=f\left(0_{V_{1}}\right)\right)$.
So $s>0$ and as such $s=m+1$ for some $m\in_{V_{1}}\omega\left(V_{1}\right)$.
Therefore $f\left(s\right)=f\left(m\right)+_{2}1_{V_{2}}$. As $s$
was minimal it holds that $m\not\in_{V_{1}}S$. and thus for $m$
\begin{equation}
\forall a\in\omega\left(V_{2}\right)\exists c\in\omega\left(V_{1}\right)\left(a<f\left(m\right)\rightarrow a=f\left(c\right)\right)
\end{equation}
 as $s\in S$ we know 
\[
\exists a\in\omega\left(V_{2}\right)\forall c\in\omega\left(V_{1}\right)\left(\left(a<_{2}f\left(s\right)\right)\wedge\left(a\neq f\left(c\right)\right)\right)
\]
let $\exists a\in\omega\left(V_{2}\right)$ be constant and receive
\begin{equation}
\forall c\in\omega\left(V_{1}\right)\left(\left(a<_{2}f\left(s\right)\right)\wedge\left(a\neq f\left(c\right)\right)\right)
\end{equation}

if $a<_{2}f\left(m\right)$ then by (1) we know that $\exists c\in\omega\left(V_{1}\right)$
that violates (2). Otherwise if $a=_{2}f\left(m\right)$ condition
(2) is violated with $c=m$ . Lastly if $a>_{2}f\left(m\right)$ it
holds that $a\geq_{2}f\left(m\right)+1_{V_{2}}=f\left(s\right)$ which
is a contradiction to $a<_{2}f\left(s\right)$ in (2). As we've received
a contradiction it must be the case that $S$ above is empty. Thus
\[
\forall b\in\omega\left(V_{1}\right)\forall a\in\omega\left(V_{2}\right)\exists c\in\omega\left(V_{1}\right)\left(a<f\left(b\right)\rightarrow a=f\left(c\right)\right)
\]

$\blacksquare$

\paragraph{Remark:}

Please note that the condition $V_{2}$ is a set of $V_{1}$ gives
us $V_{1}\underset{I}{\overset{f}{\longrightarrow}}V_{2}$. But in
fact it's a much much stronger assertion, as in this case $f$ is
also in $V_{1}$.

Specifically, one could use this fact to define induction on $\left[a\right]$
for $a\in\omega\left(V_{2}\right)\cap Im\left(f\right)$ and pull
back the argument by $f^{-1}$ to an induction on $\omega\left(V_{1}\right)$.
And, use this fact to use inductive arguments on $Im\left(f\right)$
as seen above.

\subsubsection{Overspill principle \label{subsec:Overspill-principle}}

\paragraph{Theorem:}

Let $V,\in_{V}$ be a $ZF$ model and $M,\in_{M}$ and another $ZF$
let $a_{1},\ldots,a_{n}\in M$ and $\phi$ be a formula set theory
language. s.t
\begin{itemize}
\item $M$ is a set model of $V$
\item $\omega\left(M\right)\neq\omega\left(V\right)$ 
\item $V\underset{I}{\overset{f}{\longrightarrow}}M$
\item for every $x\in\omega\left(V\right)$ it holds $M\models\phi\left(f\left(x\right),a_{1},\ldots,a_{n}\right)$
\end{itemize}
then exactly one of the following holds:
\begin{itemize}
\item $\forall n\in\omega\left(M\right)$ it holds $M\models\phi\left(n,a_{1},\ldots,a_{n}\right)$
\item exist $k\in\omega\left(M\right)$ which is non standard w.r.t $V$
s.t $\forall n<k$ it holds $M\models\phi\left(n,a_{1},\ldots,a_{n}\right)$
\end{itemize}

\paragraph{Proof:}

Within $M$ look at the set $A=\left\{ n\in\omega\left(M\right)\,|\,\neg\phi\left(n,a_{1},\ldots,a_{n}\right)\right\} $.
If the set $A$ is empty then $\forall n\in\omega\left(M\right)$
it holds $M\models\phi\left(n,a_{1},\ldots,a_{n}\right)$ and the
theorem holds true with the first condition. Otherwise $A$ isn't
empty and as a subset of natural numbers, it has a minimum. Let $k=\min A$
as for every $x\in\omega\left(V\right)$ it holds $M\models\phi\left(f\left(x\right),a_{1},\ldots,a_{n}\right)$
it follows that $k$ must be non standard w.r.t $V$. As $k$ was
the minimum of $A$ it holds $\forall n<k$ it holds $M\models\phi\left(n,a_{1},\ldots,a_{n}\right)$
and the theorem holds true with the second condition. $\blacksquare$

\subsubsection{Non-standard number definition \label{subsec:non-standard-number-definition}}

\paragraph{Definition:}

For a set model $M$ of $V$ we know by \ref{subsec:basic-properties-of-extension}
that $f:\omega\left(V\right)\rightarrow\omega\left(M\right)$ exist
s.t $V_{1}\underset{I}{\overset{f}{\longrightarrow}}V_{2}$. $j\in\omega\left(M\right)$
is called \textbf{non-standard }w.r.t $V$ if $j\not\in f\left(\omega\left(V\right)\right)$.

\paragraph{Definition:}

For a set model $M$ of $V$ we know by \ref{subsec:basic-properties-of-extension}
that $f:\omega\left(V\right)\rightarrow\omega\left(M\right)$ exist
s.t $V_{1}\underset{I}{\overset{f}{\longrightarrow}}V_{2}$. A model
$M$ of $ZFC$ in $V$ has a standard $\omega$ w.r.t $V$ if $f\left(\omega\left(V\right)\right)=\omega\left(M\right)$.
Or, in other words, $M$ has no non-standard numbers.

\subsubsection{Absoluteness of Turing machines in sub-models\label{subsec:absoluteness-of-Turing}}

\paragraph{Theorem:}

Let $V_{1},V_{2}$ be two models of $ZF$ s.t $V_{2}$ is a set of
$V_{1}$ $V_{1}\underset{I}{\overset{f}{\longrightarrow}}V_{2}$ and
let $T_{1}\in\omega\left(V_{1}\right)$ represent a coding of a Turing
machine. Let $z\in\omega\left(V_{1}\right)$ be any number. 

Denote $R\left(T_{1},z\right)\in\omega\left(V_{1}\right)$ to be the
coded state of the machine $T_{1}$ after $z$ steps (in $V_{1}$)
and $R\left(f\left(T_{1}\right),f\left(z\right)\right)$ be the coded
state of the machine $f\left(T_{1}\right)$ after $f\left(z\right)$
steps (in $V_{1}$).

Then 
\[
f\left(R\left(T_{1},z\right)\right)=R\left(f\left(T_{1}\right),f\left(z\right)\right)
\]

\paragraph{Proof (sketch):}

Denote $J$ to be the operation of running the machine in a specific
status one more step, i,e the function that 
\[
J\left(R\left(T,z\right)\right)=R\left(T,z+1\right)
\]
 for any TM $T\in\omega\left(V_{1}\right)$ and any number $z\in\omega\left(V_{1}\right)$.
The function $J$ can be expressed using arithmetic operations. The
proof is done by induction on $\omega\left(V_{1}\right)$:

The base case is $R\left(T_{1},0_{V_{`}}\right)$ as $T_{1}$ is a
machine that is coded by a text to integer coding , it holds that
$f\left(T_{1}\right)$ is a text to integer coding of the same machine
in $\omega\left(V_{2}\right)$. Thus we get 
\[
f\left(R\left(T_{1},0_{V_{1}}\right)\right)=R\left(f\left(T_{1}\right),f\left(0_{V_{1}}\right)\right)=R\left(f\left(T_{1}\right),0_{V_{2}}\right)
\]
 for $z+1$ we recall that $J$ is an arithmetic function and thus
$f\left(J\left(x\right)\right)=J\left(f\left(x\right)\right)$ for
any $x\in\omega\left(V_{1}\right)$ and thus 
\begin{align*}
f\left(R\left(T,z+1\right)\right) & =f\left(J\left(R\left(T_{1},z\right)\right)\right)\\
=f\left(J\left(R\left(T_{1},z\right)\right)\right) & =J\left(f\left(R\left(T_{1},z\right)\right)\right)\\
=J\left(R\left(f\left(T_{1}\right),f\left(z\right)\right)\right) & =R\left(f\left(T_{1}\right),f\left(z\right)+1\right)\\
 & =R\left(f\left(T_{1}\right),f\left(z+1\right)\right)
\end{align*}
Now for the induction, denote the set 
\[
S=\left\{ z\in\omega\left(V_{1}\right)\,|\,f\left(R\left(T_{1},z\right)\right)\neq R\left(f\left(T_{1}\right),f\left(z\right)\right)\right\} 
\]
 the set $S$ is a subset of natural numbers definable in $V_{1}$
(as $V_{2}$ is a set of $V_{1}$). If $S$ isn't empty, it must have
a minimum let $z'$ be that minimum. $z'\neq0_{V_{1}}$ as we've shown
that 
\[
f\left(R\left(T_{1},0_{V_{1}}\right)\right)=R\left(f\left(T_{1}\right),0_{V_{2}}\right)
\]
 if $z'=z''+1$ then $z''\not\in S$ and thus 
\[
f\left(R\left(T_{1},z''\right)\right)=R\left(f\left(T_{1}\right),f\left(z''\right)\right)
\]
 and thus we get that $z'=z''+1\not\in S$ as 
\[
f\left(R\left(T_{1},z''+1\right)\right)=R\left(f\left(T_{1}\right),f\left(z''+1\right)\right)
\]
 Therefore, $S$ above must be empty and thus 
\[
f\left(R\left(T_{1},z\right)\right)=R\left(f\left(T_{1}\right),f\left(z\right)\right)
\]
for every $z\in\omega\left(V_{1}\right)$.$\blacksquare$

\paragraph{Remark:}
\begin{enumerate}
\item Please note that $V_{2}$ may may more numbers which are not in $Im\left(f\right)$.
In this case, it may be that for a TM $T_{1}$, $T_{1}$ will not
halt in $V_{1}$ but will halt in $V_{2}$ . How every for any step
in $Im\left(f\right)$ the running in $V_{1}$ and $V_{2}$ will agree.
\item The ``sketch'' part of this proof is the fact that the coding of
$R$ and the $J$ operation wasn't fully defined. As I trust the reader
is familiar with such constructions, I don't see added value in elaborating. 
\item Please note the importance of the assumption: $V_{2}$ is a set of
$V_{1}$. The fact $V_{1}\underset{I}{\overset{f}{\longrightarrow}}V_{2}$
alone isn't enough for this proof (as we need to use induction on
$\omega\left(V_{1}\right)$). However, the above theorem still holds
true in the case of $V_{1}\underset{I}{\overset{f}{\longrightarrow}}V_{2}$
but as it won't be used it isn't shown.
\end{enumerate}

\subsubsection{Absoluteness of $\mathcal{L}\left(T\right)$ }

\paragraph{Theorem:}

Let $V_{1},V_{2}$ be two models of $ZF$ s.t $V_{2}$ is a set of
$V_{1}$, $V_{1}\underset{I}{\overset{f}{\longrightarrow}}V_{2}$
and let $T_{1}\in\omega\left(V_{1}\right)$ be a coding of a TM the
recognizes $T$ , an effective consistent set of axioms (we assume
here that $T$ is consistent according to both $\omega\left(V_{1}\right)$
and $\omega\left(V_{2}\right)$) .Let $\mathcal{L}_{1}\left(T_{1}\right),\mathcal{L}_{2}\left(f\left(T_{1}\right)\right)$
be $\mathcal{L}\left(T\right)$ computed within $V_{1},V_{2}$ respectively. 

Then $f\left(\mathcal{L}_{1}\left(T_{1}\right)\right)=\mathcal{L}_{2}\left(f\left(T_{1}\right)\right)$.

\paragraph{Intuition:}

Recall the definition of $\mathcal{L}\left(T\right)$ given in \ref{subsec:Chaitin's-incompleteness-theorem}:
$\mathcal{L}\left(T\right)$ to be the minimal number s.t $L>\log\left(L\right)+C$.
Therefore, as long as $C$ is interpreted the same in both models,
$\mathcal{L}\left(T\right)$ will also be the same. 

Moreover, $C$ contains:
\begin{enumerate}
\item Representation of a TM which computes $T$ (which is assumed to be
the same).
\item A representation of the machine $M$ which goes over the proofs from
$T$ and finds the first proof of the claim in $\left\{ L\leq Kol\left(x\right)|x\in\omega\left(V\right)\right\} $
(for a given fixed $L$).
\end{enumerate}
As both of these are absolute in sub-models (as seen in \ref{subsec:absoluteness-of-Turing})
$C$ must be interpreted the same.

\paragraph{Proof (sketch):}

The reader may want to think of the case $\omega\left(V_{1}\right)=\mathbb{N}$
at first. The idea of this proof is that all of \ref{subsec:Chaitin's-incompleteness-theorem}
construction could be done in a bounded set of $\omega\left(V_{1}\right)$
and $\omega\left(V_{2}\right)$ behaves ``the same'' within on such
subsets (lower parts). Here are some details:

Recall the proof of \ref{subsec:Chaitin's-incompleteness-theorem},
the proof creates a Turing machine $M$ which goes over over all elements
of $\omega\left(V_{1}\right)$ until it finds the first element $w$
which proves a claim in the set $\left\{ L\leq Kol\left(x\right)|x\in\omega\left(V\right)\right\} $.
\begin{itemize}
\item Proofs are a list of statements in first order logic, each statement
can be either an axiom or derived from the previous statements.
\item Such a proof can be coded into an integer be text to integer coding.
Such a coding require only the operations of addition, multiplication,
exponentiation.
\item The proof could be verified by a Turing machine, which runs the same
(by \ref{subsec:absoluteness-of-Turing}) in both $V_{1}$ and $V_{2}$.
\item Therefore the proofs interprets the same in both $V_{1}$ and $V_{2}$.
\begin{itemize}
\item i.e if $T\underset{p}{\vDash}\phi$ in $V_{1}$ then $f\left(T\right)\underset{f\left(p\right)}{\vDash}f\left(\phi\right)$
in $V_{2}$. where $T,p,\phi$ are coding of a theory, proof and a
statement respectively.
\end{itemize}
\item Let $w\in\omega\left(V_{1}\right)$ be the first proof of a statement
in $\left\{ L\leq Kol\left(x\right)|x\in\omega\left(V_{1}\right)\right\} $
in $V_{1}$. The proof $f\left(w\right)$ is also a proof of a statement
in $\left\{ L\leq Kol\left(x\right)|x\in\omega\left(V_{2}\right)\right\} $
in $V_{2}$ and due to the one to one correspondence of $f$, $f\left(w\right)$
must be the first such proof in $\omega\left(V_{2}\right)$ as well.
\item Therefore the operation of machine $M$ will work ``the same'' and
return $x$ in $V_{1}$ and $f\left(x\right)$ in $V_{2}$.
\item for $L_{1}\in\omega\left(V_{1}\right)$,$L_{2}\in\omega\left(V_{2}\right)$
we've created two Turing machines $M_{1}$ in $V_{1}$ and $f\left(M_{1}\right)$
in $V_{2}$ that computes $x',f\left(x'\right)$ in $V_{1},V_{2}$
respectively. Therefore in both $V_{1},V_{2}$ the two inequalities
must hold:
\[
L_{1}\leq\log\left(L_{1}\right)+C
\]
\[
L_{2}\leq\log\left(L_{2}\right)+f\left(C\right)
\]
 and this the minimum number that violates them must be ``the same''
i.e 
\[
f\left(\mathcal{L}_{1}\left(T_{1}\right)\right)=\mathcal{L}_{2}\left(f\left(T_{1}\right)\right)
\]
\end{itemize}

\subsubsection{Reduction of machine number \label{subsec:reduction-of-machine}}

\paragraph{Theorem:}

Let $T$ be a effective consistent theory containing the axioms of
$ZF$ ($T$ may contain some more consistencies axioms as well). Let
$C$ be the size of the universal TM.

Assume
\begin{itemize}
\item $M$ is a model of $T\cup con-\left(T\right)$
\item $\phi\left(n,k\right)$ is a two variable formula that defines a function
in $ZF$ and hold the propriety condition
\item $t\in\omega\left(M\right)$
\item $x\in\omega\left(M\right)$ codes a TM which computes the mapping
$n\rightarrow k$ s.t $\phi\left(n,k\right)$ holds for every $n<t$
in $\left\lceil \log\left(n+1\right)+\log\left(k+1\right)+1\right\rceil $
computational steps.
\item and let $\overset{\,o}{Z}\in\omega\left(M\right)$ s.t $\overset{\,o}{Z}>\mathcal{L}\left(T\right)+C$
.
\end{itemize}
Then, a model $N$ of $T$ exist s.t:
\begin{itemize}
\item $N$ is a set of $M$
\item $M\underset{I}{\overset{f}{\longrightarrow}}N$ for a function $f:\omega\left(M\right)\rightarrow\omega\left(N\right)$ 
\item $\exists y,c\in\omega\left(N\right)$ s.t $y\leq f\left(\overset{\,o}{Z}\right)$
and $y$ codes a TM which computes for every $n<f\left(t\right)$
the value of $k$ s,t $\phi\left(n,k\right)$ in $\left\lceil \log\left(n+1\right)+\log\left(k+1\right)+c\right\rceil $
steps. 
\end{itemize}

\paragraph{Proof:}

Let $M,\phi,t,x,\overset{\,o}{Z}$ be as in the theorem. If $x<\overset{\,o}{Z}$
then $M=N$ and $y=x$ and $c=1$ holds the conclusions of the theorem.
Otherwise assume $x\geq\overset{\,o}{Z}$. As $\overset{\,o}{Z}>\mathcal{L}\left(T\right)$
and $x\geq\overset{\,o}{Z}$ we get $x>\mathcal{L}\left(T\right)$.As
$M$ is a model of $T\cup con-\left(T\right)$ and by Chaitin's incompleteness
theorem\ref{subsec:Chaitin's-incompleteness-theorem} (on $M$) we
know that the statement ``$Kol\left(x\right)>\mathcal{L}\left(T\right)$''
can't be proven from $T$. For that reason the set of axioms $T\cup\left\{ "Kol\left(x\right)\leq\mathcal{L}\left(T\right)"\right\} $
is a consistent set of axioms (within $M$) thus $M$ has a model
of $T\cup\left\{ "Kol\left(x\right)\leq\mathcal{L}\left(T\right)"\right\} $.
Let $N$ be this model . As $N$ is a set model of $M$ by \ref{subsec:basic-properties-of-extension}
we know $M\underset{I}{\overset{f}{\longrightarrow}}N$. As $N$ holds
$"Kol\left(x\right)\leq\mathcal{L}\left(T\right)"$ In $\omega\left(N\right)$
exist $y'\in_{N}\omega\left(N\right)$ s.t $y'\leq f\left(\overset{\,o}{Z}\right)$
and $y'$ codes a TM that calculates $x$. The TM $y$ receives $n<f\left(t\right)$
first calculates $x$ using $y'$ and then execute $x$ on $n$. The
size of $y$ is at most $y'$ and the size of a universal TM on it's
output and hence $y<f\left(\overset{\,o}{Z}\right)$. The running
time of $y$ is the same as $x$ up to a constant hence for every
$n<f\left(t\right)$ $y$ calculates the value of $k$ s.t $\phi\left(n,k\right)$
in $\left\lceil \log\left(n+1\right)+\log\left(k+1\right)+c\right\rceil $
steps for some constant $c$. $\blacksquare$

\section{What is knowledge?\label{sec:What-is-knowledge?}}

In this section we break the sequence of the construction in order
to discuss the implications of theorem \ref{subsec:reduction-of-machine}
and how it leads to a definition of knowledge \footnote{The reader who wishes to skip this section may jump to section \ref{sec:Extension and Lattices. Second part.}
which is a sequel to section \ref{sec:Lattices}.}. 

Let's look at theorems \ref{subsec:reduction-of-machine} and \ref{subsec:binary-tree-construction}.
Together they state that given an arbitrary sequence of natural numbers
of arbitrary length, one can construct a model in which the sequence
is computable in linear time using a machine of bounded size. It is
of crucial importance to emphasize that the size of the machine bound
is \textbf{independent of the length and the numbers of the chosen
sequence}. It is philosophically unacceptable that by pure coincidence
it just ``happens'' that we can compute the chosen sequence using
a machine of a small size for every choice on the sequence. The construction
of theorem \ref{subsec:binary-tree-construction} alone gave us a
machine which represented a table of values and depended on the length
and the numbers on the sequence. Therefore there is no cognitive dissonance
when we conceive it as having full knowledge on the sequence. This
way of looking at things is incompatible with \ref{subsec:reduction-of-machine}:
it is inconceivable that the machine number still holds full knowledge
on the sequence while the size of the machine is bounded. We must
conclude that the knowledge on the sequence got transferred to the
structure of $\omega$ of the new model. So we must conclude that
the new model in \ref{subsec:reduction-of-machine} has ``learned''
the knowledge hidden in the sequence.

The following question suggests itself naturally: can we make a construction
similar to the one in theorem \ref{subsec:reduction-of-machine} for
an infinite sequence? 

Let's assume that we've successfully done this. An infinite sequence
is a function from $\omega$ to itself. Let's further assume that
this function is definable (see \ref{subsec:Definable-function})
as it must have an interpretation in different models. Let's also
assume that the definition satisfies the propriety condition (see
\ref{propriety-condition-definition}) because the interpretation
of the function in our construction is the correct interpretation
of the function in the new model. 

So, given a two variable formula $\phi\left(n,k\right)$ in the language
if set theory which defines a function in $ZF$ and which satisfies
the propriety condition and given an $x\in\omega$, we take the sequence
$\left(a_{n}\right)_{n=0}^{x}$ s.t 
\[
\forall n\,\,\,0\leq n\leq x\rightarrow\phi\left(n,a_{n}\right)
\]
 and use \ref{subsec:reduction-of-machine} to create a model $M_{x}$
in which the function is ``known'' up to $x$. Next we must ``tie''
or ``combine'' all these models together in order to create an all
encompassing model $\stackrel[x]{}{\oplus}M_{x}$. This object isn't
defined yet, but we wish the model to be such as to have a constant
$\overset{\,o}{Z}$ which bounds all TMs that compute an element in
one of the sequences. We also want every $n$ in $\omega$ to be contained
in at least one sequence. Formally we define:
\begin{equation}
\begin{matrix}\exists\overset{\,o}{Z}\in\omega\,\,\forall n\in\omega\,\,\exists x,k\in\omega\\
\left(\begin{matrix}\phi\left(n,k\right)\,\wedge\,Run\left(x,n\right)=k\,\,\wedge\,x<\overset{\,o}{Z}\end{matrix}\right)
\end{matrix}\label{eq:43}
\end{equation}
where $Run\left(x,n\right)$ denotes the function that returns the
value returned by the TM numbered $x$ on input $n$.

Note that the critical demand of running time is omitted from defintion
(\ref{eq:43}). Recall that the TMs in \ref{subsec:reduction-of-machine}
worked using linear time computation on ``legal'' inputs (i.e halted
and gave the right answer). We won't make any demands regarding running
time (or even halting) on other inputs. So, a coding $x$ of a TM
is of the right running time if 
\[
\begin{matrix}\left(\exists c\in\omega\right)\left(\forall n',k'\in\omega\right)\\
\left(\begin{matrix}\left(\phi\left(n',k'\right)\wedge Run\left(x,n'\right)=k'\right)\rightarrow\\
\left(Time\left(x,n'\right)=\left\lceil \log\left(n'+1\right)+\log\left(k'+1\right)+c\right\rceil \right)
\end{matrix}\right)
\end{matrix}
\]
 Therefore, the full definition of ``knowledge on the function $\phi$''
will be 
\begin{equation}
\begin{matrix}\exists\overset{\,o}{Z}\in\omega\,\,\forall n\in\omega\,\,\exists x,k,c\in\omega\\
\left(\begin{matrix}\phi\left(n,k\right)\,\wedge\,Run\left(x,n\right)=k\,\,\wedge\,x<\overset{\,o}{Z}\,\,\wedge\\
\left(\forall n',k'\in\omega\right)\left(\begin{matrix}\left(\phi\left(n',k'\right)\wedge Run\left(x,n'\right)=k'\right)\rightarrow\\
\left(Time\left(x,n'\right)=\left\lceil \log\left(n'+1\right)+\log\left(k'+1\right)+c\right\rceil \right)
\end{matrix}\right)
\end{matrix}\right)
\end{matrix}\label{eq:44}
\end{equation}
 where $Run\left(x,n\right)$ denotes the function that returns the
value returned by the TM numbered $x$ on input $n$ and $Time\left(x,n\right)$
returns the number of steps done in the calculation of $x$ on input
$n$. As this definition is inspired by an extrapolasion of \ref{subsec:reduction-of-machine}
we expect that definition (\ref{eq:44}) will be held by at least
some models of $ZFC$. 

Please note: 
\begin{itemize}
\item The definition of knowledge is a computation definition which is \uline{different}
from the more common definition of ``computational decision''.
\begin{itemize}
\item every function with a bounded image is clearly known, even these function
which aren't computable \textbackslash{} aren't computable in linear
time \textbackslash{} aren't computable in efficient time (in whichever
definition of efficient we may use).
\end{itemize}
\item Unlike in the traditional definition, the running time comes ``baked
in'' this definition and must be always linear.
\begin{itemize}
\item In traditional definition, the larger the running time bound the more
languages one can decide using such a running time. We don't expect
the same to be the case in knowledge as we expect linear running time
to be enough.
\end{itemize}
\item Unlike in the traditional definition, the running time must be \textbf{equal}
to linear and not just linearly bounded. 
\begin{itemize}
\item As the running time was dictated from the running time of \ref{subsec:binary-tree-construction}
with an addition of a uniform constant in \ref{subsec:reduction-of-machine}.
A running time of equal or less than linear will introduce other TMs
that work in a different fashions (which isn't our intent).
\end{itemize}
\end{itemize}
Given such $\phi$ the question of whether or not $\phi$ is known
in a model (or even if there is a model where $\phi$ is known) is
a \textbf{\uline{percolation}} conjecture. In \ref{sec:Percolation-theorems}
we will show that in a model with a worldly cardinal $\phi$ is known
for any $\phi$ two variable formula $\phi\left(n,k\right)$ (in the
language if set theory) which defines a function in $ZF$ and which
holds the propriety condition.

\section{Extension and Lattices. Second part.\label{sec:Extension and Lattices. Second part.}}

\subsection{Base model definitions\label{subsec:Base-model-definitions}}

Now we will start to define the models in question. We start with
$W$, a model of $ZFC$ with a worldly cardinal. Within $W$ exist
a countable transitive model of $ZFC+\left(\omega-con\left(ZFC\right)\right)$
called $V$. Within $V$ exist a countable (non-transitive) model
$M_{1}$ s.t $j\in\omega\left(M_{1}\right)$ exist where $j$ isn't
standard and $M_{1}$ is a model of $ZFC+j-con\left(ZFC\right)$.
These will be out base models and base on them we will use forcing
in the next section. 

\paragraph{Notation:}

Let $W$ be a model of $ZFC$ with a worldly cardinal. 

\paragraph{Theorem:}

$W$ holds $\omega-con\left(ZFC\right)$ 

\paragraph{Proof:}

Denote for $k\in ON\left(W\right)$ denote $V_{k}^{W}$ to be the
Von-neumann universe of $W$ and let $k'\in ON\left(W\right)$ be
worldly (i.e $V_{k}^{W},\in_{W}$ is a ZFC model). We will prove $\forall t\in\omega\left(W\right)\,\,W\models t-con\left(ZFC\right)$
by induction over $t$. 
\begin{itemize}
\item Base case $t=1$. As $V_{k}^{W},\in_{W}$ is a ZFC model $W$ has
a set model of $ZFC$ and as such can't prove a contradiction from
$ZFC$ and thus in $W$, $ZFC$ is consistent and so $W$ holds $1-con\left(ZFC\right)$.
\item Assume that $W$ holds $t-con\left(ZFC\right)$. As the property $t-con\left(ZFC\right)$
can be expressed as a property of natural numbers (the set of axioms
can't prove a contradiction) and as $\omega\left(W\right)=\omega\left(V_{k}^{W}\right)$
we get that $V_{k}^{W}$ as a set model also hold $t-con\left(ZFC\right)$.
As such $V_{k}^{W}$ has a a sequence of $ZFC$ models $M'_{1},M'_{2},M'_{3},...,M'_{t}$
within $V_{k}^{W}$ s.t $M'_{1}$ is a set of $V_{k}^{W}$ and each
model is a set of its previous. Thus $W$ has a sequence $V_{k}^{W},M'_{1},M'_{2},M'_{3},...,M'_{t}$
of $ZFC$ models s.t $V_{k}^{W}$ is a set of $W$ and every model
is a set of its previous. Thus $W\models\left(t+1\right)-con\left(ZFC\right)$ 
\end{itemize}
As we got that for all $t\in\omega\left(W\right)$ $W\models t-con\left(ZFC\right)$
we know by definition that $W$ holds $\omega-con\left(ZFC\right)$
.$\blacksquare$

\paragraph{Corollary:}

As $W$ holds $\omega-con\left(ZFC\right)$ and as $\omega-con\left(ZFC\right)$
is a property that can be expressed as a property of natural numbers
we get from the same argument that $V_{k}^{W}$ also holds $\omega-con\left(ZFC\right)$.

\subsubsection{$V$ construction }

\paragraph{Theorem:}

$W$ has a countable transitive set model of $ZFC+\left(\omega-con\left(ZFC\right)\right)$.

\paragraph{Proof:}

As seen in previously in \ref{subsec:Base-model-definitions} $ZFC+\left(\omega-con\left(ZFC\right)\right)$
axioms are consistent in $W$ and have a model $V_{k}^{W}$ of them.
As the language of set theory is countable by Lowenheim Skolem theorem
\ref{subsec:Lowenheim-Skolem-theorem} exists $X$ a countable elementary
sub model.

$X$, as a subset of $V_{k}^{W}$, is also well founded w.r.t $\in_{W}$
. $X$ is countable but might not be transitive. By Mostowski\textquoteright s
collapsing theorem \ref{subsec:Mostowski=002019s-Collapsing-Theorem},
we know that $X$ can be collapsed to a transitive set $V$ and so
$V,\in_{W}$ is a countable transitive set model of $ZFC+\left(\omega-con\left(ZFC\right)\right)$.
$\blacksquare$ 

\paragraph{Notation:}

Let $V$ be the model a countable transitive set model $ZFC+\left(\omega-con\left(ZFC\right)\right)$
within $W$.

\paragraph{Lemma:}

It holds that $\omega\left(V\right)=\omega\left(W\right)$.

\paragraph{Proof:}

As $V$ is a set model of $W$ we know by \ref{subsec:basic-properties-of-extension}
that $W\underset{I}{\overset{f}{\longrightarrow}}V$ and as $\emptyset$
and the successor operation interprets the same in $W$ and $V$ we
know that $\omega\left(W\right)\subseteq\omega\left(V\right)$. As
$V$ is transitive, if $\omega\left(W\right)\subsetneq\omega\left(V\right)$
then in $W$ the set $\omega\left(W\right)\smallsetminus\omega\left(V\right)$
must contain an infinite decreasing sequence as for every number $x\in_{W}\omega\left(W\right)\smallsetminus\omega\left(V\right)$
the number $x-k$ for $k\in\omega\left(W\right)$ is also in $\omega\left(W\right)\smallsetminus\omega\left(V\right)$.
As $W$ is a $ZFC$ model it can't contain an infinite decreeing sequence
and thus $\omega\left(V\right)=\omega\left(W\right)$. $\blacksquare$

\subsubsection{$M_{1}$ construction \label{subsec:M_1-construction}}

\paragraph{Theorem:}

$V$ has a countable model $M_{1}$, $j\in\omega\left(M_{1}\right)$
non standard w.r.t $V$ and $M_{1}$ is a model of $ZFC+\left(j-con\left(ZFC\right)\right)$
but not of $ZFC+\left(\left(j+1\right)-con\left(ZFC\right)\right)$

\paragraph{Proof:}

Extend the language of set theory to include one more variable $j'$.
Build the following set of axioms over the extended language:
\begin{enumerate}
\item all axioms of $ZFC$
\item for every number $n\in\omega\left(V\right)$ add the following axioms:
\begin{enumerate}
\item $j'>n$ 
\item $\left(j'-con\left(ZFC\right)\right)$
\end{enumerate}
\item $\exists n'\in\omega\,\neg\left(n'-con\left(ZFC\right)\right)$
\end{enumerate}
As every finite set of these axioms is consistent (with a choice of
a large enough $j'$ from $\omega\left(V\right)$) we know by the
compactness theorem \ref{subsec:Compactness-theorem} that the whole
set of axioms is consistent. As the axioms set is consistent let $M_{1}'$
be a model. As the extended language is countable, by Lowenheim Skolem
theorem \ref{subsec:Lowenheim-Skolem-theorem} we know that a countable
model $M_{1}$ exist in which $j'\in\omega\left(M_{1}\right)$ exist
s.t for every $n\in\omega\left(V\right)$ $j'>n$ and $M_{1}$ holds
$\left(j'-con\left(ZFC\right)\right)$. As the axiom $\exists n'\in\omega\,\neg\left(n'-con\left(ZFC\right)\right)$
holds we now that $M_{1}$ doesn't hold $\omega-con\left(ZFC\right)$
denote $j$ to be the maximal $j'$ s.t $\left(j'-con\left(ZFC\right)\right)$
holds in $M_{1}$$\blacksquare$

\paragraph{Notation:}

Denote the above model $M_{1}$. Denote the collection of $ZFC$ model
within $M_{1}$ to be $Models_{M_{1}}\left(ZFC\right)$.

\paragraph{Corollary:}

As $M_{1}$ is a countable set of $V$ and as being a model is a property
(which some sets in $M_{1}$ have and some don't) it following that
$Models_{M_{1}}\left(ZFC\right)$ is a countable set in $V$.

\subsection{Forcing POS over $V$ \label{subsec:Forcing-POS-over-V}}

Now given $W,V,M_{1}$ we define the forcing conditions and generic
filters over $V$. $W$ will be used as a would model on which we
will build our generic filters and $V$ will be the model being extended.
$M_{1}$ will be a prat of the conditions.

\paragraph{Notation:}

For a $ZFC$ model $V'$ denote the power set of a set $A$ in $V'$
to be $POW_{V}\left(A\right)$

\paragraph{Definition:}

For a $ZFC$ model $V'$ a finite ordered set is a function $f$ with
domain $\left[a\right]=\left\{ i\in\omega\left(V'\right)|i\leq a\right\} $
where $a\in\omega\left(V'\right)$.

\paragraph{Definition:}

Define the following POS within $V$:
\[
P_{0}=\left\{ \left(f_{1},f_{2},f_{3}\right)\,|\,\begin{matrix}f_{1}:Models_{M_{1}}\left(ZFC\right)\rightharpoonup Pow_{V}\left(\omega\left(V\right)\right)\\
f_{2}:Models_{M_{1}}\left(ZFC\right)\rightharpoonup Pow_{V}\left(\omega\left(V\right)\right)\\
f_{2}:Models_{M_{1}}\left(ZFC\right)\rightharpoonup Pow_{V}\left(\omega\left(V\right)\right)\\
\text{\ensuremath{f_{1},f_{2},f_{3}} are partial functions with a finite (by \ensuremath{V}) domain.}\\
Dom\left(f_{1}\right)=Dom\left(f_{2}\right)=Dom\left(f_{3}\right).\\
\text{for \ensuremath{M\in\ensuremath{Models_{M_{1}}\left(ZFC\right)}} if \ensuremath{M\in Dom\left(f_{1}\right)} then}\\
f_{1}\left(M\right)\text{ is a finite (by \ensuremath{V}) ordered set of elements within \ensuremath{M} i.e \ensuremath{\left\{ e\in M|e\in f_{1}\left(M\right)\right\} }}\\
f_{2}\text{\ensuremath{\left(M\right)} is a finite (by \ensuremath{V}) ordered set of of formulas in the language }\\
\text{of set theory with varibles }\ensuremath{\left\{ e_{i}\right\} _{i\in f_{1}\left(M\right)}}\\
f_{3}\left(M\right)\text{ is an element of \ensuremath{\omega\left(M\right)} }
\end{matrix}\right\} 
\]
along with the partial order $f\leq_{P_{0}}g$ if $Dom\left(f_{1}\right)\subset Dom\left(g_{1}\right)$
and $g_{1}\restriction_{Dom\left(f\right)}=f_{1}$ and $g_{2}\restriction_{Dom\left(f\right)}=f_{2}$
and $g_{3}\restriction_{Dom\left(f\right)}=f_{3}$.

\paragraph{Explanation:}

$P_{0}$ is a set of partial functions. A function $f$ is a function
from a finite subset of $Models_{M_{1}}\left(ZFC\right)$ that returns
for each $M$ in its domain:
\begin{itemize}
\item a finite (by $V$) set of elements of $M$.
\item a finite (by $V$) set of formulas in the language of set theory with
free variables corresponding to the above set.
\item as, by text to integer coding, every formula can by coded by an integer
and thus the function is into $\omega\left(V\right)$.
\begin{itemize}
\item please note that the formulas must by in $V$.
\item please note that there is no consistency requirement of these formulas 
\end{itemize}
\end{itemize}

\paragraph{Lemma:}
\begin{enumerate}
\item Within $W$ a filter $G_{0}\subset P_{0}$ exist which is a generic
filter over $V$
\item $V\left[G_{0}\right]$ is a countable transitive model of $ZFC+\left(\omega-con\left(ZFC\right)\right)$.
\end{enumerate}

\paragraph{Proof:}

By application of theorem 1 in \ref{subsec:Forcing} we get the filter
$G_{0}$ and by applying theorem 4 we get that $V\left[G_{0}\right]$
is a countable transitive model of $ZFC$. As $V\left[G_{0}\right]$
is transitive $\omega\left(V\left[G_{0}\right]\right)=\omega\left(V\right)$
and as $\left(\omega-con\left(ZFC\right)\right)$ is a set of formulas
that can be expressed as natural numbers and as $V$ holds $ZFC+\left(\omega-con\left(ZFC\right)\right)$
we get that $V\left[G_{0}\right]$ is a countable transitive model
of $ZFC+\left(\omega-con\left(ZFC\right)\right)$. $\blacksquare$

\paragraph{Notation remark:}

For a model $A,\in_{A}$ in $M_{1}$ we will use the notation $G_{0}\left(A\right)$
to symbolize the value of $G_{0}$ as a function on $A$. As $G_{0}$
is a collection of function (and not just a single function) it may
be unclear. However, as every two function in $G_{0}$ that have $A$
in their domain must agree on their assigned value on $A$ we can
view $\cup G_{0}$ as one big function that gives the value $G_{0}\left(A\right)$. 

\paragraph{Remark:}

Given a model $A$ in $M_{1}$ , $G_{0}\left(A\right)$ is defined.
As the set $E$ of parital function in $P_{0}$ which are defined
on $A$ is a dense set in $V$ and as such $G_{0}$ as a generic filter
must intersect it. As such , $G_{0}\left(A\right)$ is defined.

\subsubsection{Formula reduction definition\label{subsec:Formula-reduction-definition}}

\paragraph{Definition:}

Let $M$ be a model of $ZFC$ which is a set of $V$. For a set of
formulas $A$ in $V$ (which may have free variables in them) and
$k\in\omega\left(V\right)$
\begin{itemize}
\item we say that $A$ is $k$ consistent in $M$ from $ZFC$ if , $M$'s
arithmetic holds
\[
k-con\left(ZFC\cup A\right)
\]
In other words, $M$ has a sequence of $k$ models, each model is
a set of its previous and all hold $ZFC\cup A$. 
\end{itemize}

Please note that unlike $\left\{ n-con\left(ZFC\right)\,|\,n\in\omega\left(V\right)\right\} $the
set $ZFC\cup n-con\left(ZFC\right)$ for $n\in\omega\left(V\right)$
can be defined in $M$. 

\paragraph{Definition:}

Let a model $M$ of $ZFC\cup\left\{ n-con\left(ZFC\right)\,|\,n\in\omega\left(V\right)\right\} $,
let $a,b,c\in\omega\left(M\right)$ and let a set of constant $\left\{ e_{i}\right\} _{i\in\left[c\right]}$
and two function $f_{1}:\left[a\right]\rightarrow\omega\left(M\right)$
and $f_{2}:\left[b\right]\rightarrow\omega\left(M\right)$ which represent
a text to integer coding of formulas in $M$ with the constants in
$\left\{ e_{i}\right\} _{i\in\left[c\right]}$. Assume that the formulas
in $f_{1}$ are $k$ consistent in $M$ from $ZFC$ for every $k\in\omega\left(V\right)$
(i.e 
\[
k-con\left(ZFC\cup f\left(\left[a\right]\right)\right)
\]
 in $M$ for every $k\in\omega\left(V\right)$). Define the reduction
function by recursion over $b$ by the following recursive algorithm:
\begin{enumerate}
\item If $f_{2}$ isn't an empty function, look at the formula $f_{2}\left(0\right)$:
\begin{enumerate}
\item If 
\[
k-con\left(ZFC\cup f\left(\left[a\right]\right)\cup f_{2}\left(0\right)\right)
\]
 in $M$ for every $k\in\omega\left(V\right)$ define $f_{1}^{new}=f_{1}^{old}\cup\left(Dom\left(f_{1}\right),f_{2}\left(0\right)\right)$
and for $0\leq n\leq b-1$ define $f_{2}^{new}\left(n\right)=f_{2}^{old}\left(n+1\right)$.
Go to back to step $\left(1\right)$ with $\left(f_{1}^{new},f_{2}^{new}\right)$.
\item Otherwise, define $f_{1}^{new}=f_{1}^{old}\cup\left(Dom\left(f_{1}\right),\neg f_{2}\left(0\right)\right)$
and for $0\leq n\leq b-1$ define $f_{2}^{new}\left(n\right)=f_{2}^{old}\left(n+1\right)$.
Go to back to step $\left(1\right)$ with $\left(f_{1}^{new},f_{2}^{new}\right)$.
\end{enumerate}
\item If $f_{2}$ is the empty function define $Reduce_{M}\left(f_{1},f_{2}\right)=f_{1}$.
\end{enumerate}

\paragraph{Explanation:}

The process of $Reduce$ takes two lists of axioms. Where the first
list $f_{1}$ is assumed to be consistent. It adds axioms from $f_{2}$,
one at a time by their order. As axioms can be added as long as it
doesn't create a contradiction. If the axiom does create a contradiction
the negation of the axiom is added and the process goes on to the
next axiom in the list of $f_{2}$. By the end of this process (as
the lists are finite) we get a consistent list of axioms called $Reduce_{M}\left(f_{1},f_{2}\right)$
and every axiom in $f_{2}$ is either listed in $Reduce_{M}\left(f_{1},f_{2}\right)$
or its negation is listed there.

\subsection{Lattice construction\label{subsec:Lattice-construction}}

In this part, within $V\left[G_{0}\right]$, we will now construct
the set of models $\left\{ M_{i}\right\} _{i\in\omega\left(V\right)}$
on which we will later define a limit model. Please recall that $M_{1}$
is already constructed along with $j\in\omega\left(M\right)$ which
is non standard w.r.t $V$ and $M_{1}$ holds $j-con\left(ZFC\right)$.
As $M_{1}$ is a countable model (and therefore countable set), fix
a numbering (in $V$) of $M_{1}$.

\subsubsection{$M_{i}$'s construction}

The definition of $M_{i}$ will be by induction. For $i\in\omega\left(V\right)$
s.t $i\geq1$ we assume we've defined $M_{i},var_{i},Formulas_{i}$
s.t
\begin{itemize}
\item $M_{i}$ is a $ZFC$ model.
\item $M_{i}$ is a set model of $V$.
\item $M_{i}$ holds $ZFC\cup\left\{ n-con\left(ZFC\right)\,|\,n\in\omega\left(V\right)\right\} $.
\item $M_{i}$ holds $ZFC\cup Formulas_{i}\cup\left\{ n-con\left(ZFC\cup Formulas_{i}\right)\,|\,n\in\omega\left(V\right)\right\} $.
\item for $i>1$ $M_{i}$ is an element $M_{1}$
\item $var_{i}$ is a finite (by $V$) set of variables of set in $M_{i}$.
\begin{itemize}
\item For the case $i=1$ define $var_{i}$ to be the empty set.
\item As $var_{i}$ is a finite set (by $V$) of elements in $M_{i}$, $var_{i}$
is also a set of $M_{i}$.
\end{itemize}
\item $Formulas_{i}$ is a finite (by $V$) set of formulas (in $V$) with
variables in $var_{i}$.
\begin{itemize}
\item For the case $i=1$ define $Formulas_{i}$ to be an empty set of formulas.
\item As $Formulas_{i}$ is a finite set (by $V$) of elements in $M_{i}$
and each formula is a number in $\omega\left(V\right)$, $Formulas_{i}$
is also a set of $M_{i}$.
\end{itemize}
\end{itemize}
From the generic filter $G_{0}$ we receive:
\begin{itemize}
\item $\left(G_{0}\left(M_{i}\right)\right)_{1}$ is a finite (by $V$)
set of variables of set in $M_{i}$.
\item $\left(G_{0}\left(M_{i}\right)\right)_{2}$ is a finite set of formulas
(in $V$) with variables in $\left(G_{0}\left(M_{i}\right)\right)_{1}$.
Each formula is a formula in $V$.
\item $\left(G_{0}\left(M_{i}\right)\right)_{3}$ is a number in $\omega\left(M_{i}\right)$.
\end{itemize}

\paragraph{Definitions:}

Given $i+1$ define the following:

Denote 
\begin{itemize}
\item $var_{i+1}=var_{i}\cup\left(G_{0}\left(M_{i}\right)\right)_{1}$.
\item $Formulas_{i+1}=Reduce_{M_{i}}\left(Reduce_{M_{i}}\left(\emptyset,Formulas_{i}\right),\left(G_{0}\left(M_{i}\right)\right)_{2}\right)$.

Where $Reduce_{M_{i}}$ is the reduction function defined in \ref{subsec:Formula-reduction-definition}.
\end{itemize}

\paragraph{Axioms of $M_{i+1}$:}

We define the extended language to be the language of set theory along
with constants for each variable in $var_{i+1}$. The following list
of axioms is 

We demand the following axioms:
\begin{enumerate}
\item $ZFC$
\item All formulas in $Formulas_{i+1}$ must hold (as formulas with the
appropriate constant in $var_{i+1}$). and 
\item for every $j\in\omega\left(V\right)$ 
\[
j-con\left(ZFC\cup Formulas_{i+1}\right)
\]
 must hold.
\begin{enumerate}
\item Denote $k$ to be the maximal number s.t 
\[
k-con\left(ZFC\cup Formulas_{i+1}\right)
\]
 is consistent in $M_{i}$. Such $k$ must be non standard w.r.t $V$.
\end{enumerate}
\end{enumerate}

\paragraph{Remark:}
\begin{enumerate}
\item As $M_{i}$ is consistent with 
\[
ZFC\cup Formulas_{i}\cup\left\{ n-con\left(ZFC\cup Formulas_{i}\right)\,|\,n\in\omega\left(V\right)\right\} 
\]
 and as $Formulas_{i+1}$ was constructed to be consistent with 
\[
ZFC\cup Formulas_{i+1}\cup\left\{ n-con\left(ZFC\cup Formulas_{i+1}\right)\,|\,n\in\omega\left(V\right)\right\} 
\]
 in $M_{i}$ we know that if $k$ is the maximal number s.t 
\[
k-con\left(ZFC\cup Formulas_{i+1}\right)
\]
 is consistent in $M_{i}$, such a $k$ must m be non standard w.r.t
$V$ by the overspill principle \ref{subsec:Overspill-principle}
and such $k$ must exist as $M_{1}$ didn't hold $\omega-con\left(ZFC\right)$
(and consequently all $M_{i}$ won't hold $\omega-con\left(ZFC\right)$
).
\item The demand every $j\in\omega\left(V\right)$ 
\[
j-con\left(ZFC\cup Formulas_{i+1}\right)
\]
 must hold can't be stated inside $M_{i}$ (as $M_{i}$ doesn't have
access to $\omega\left(V\right)$) but given that $k$ of axioms 3
is non standard w.r.t $V$ (which it is by our construction) we can
demand 
\[
k-con\left(ZFC\cup Formulas_{i+1}\right)
\]
 and the demand that for every $j\in\omega\left(V\right)$ 
\[
j-con\left(ZFC\cup Formulas_{i+1}\right)
\]
 follows that statement. 
\end{enumerate}

\paragraph{Lemma:}

The axioms 1+2+3 above can be demanded in $M_{i}$ and are consistent
in it.

\paragraph{Proof:}

Let $k$ be the maximal number s.t 
\[
k-con\left(ZFC\cup Formulas_{i+1}\right)
\]
 is consistent in $M_{i}$.

The model $M_{i}$ is assumed to be consistent with 
\[
j-con\left(ZFC\cup Formulas_{i}\right)
\]
for every $j\in\omega\left(V\right)$ . $Formulas_{i+1}$ is a $j$
consistent set of axioms for every $j\in\omega\left(V\right)$ (as
$Formulas_{i+1}$ was chosen such). Thus, for every $j\in\omega\left(V\right)$
\[
j-con\left(ZFC\cup Formulas_{i+1}\right)
\]
is consistent in $M_{i}$. Such a $k$ therefore must be non-standard
w.r.t $V$. 

By previous remark in order to show that we can demand axiom 3 in
$M_{i}$ (as $\omega\left(V\right)$ can't be defined within $M_{i}$)
it is suffice to show that 
\[
k-con\left(ZFC\cup Formulas_{i+1}\right)
\]
 is consistent in $M_{i}$ and that $k$ is non standard w.r.t $V$.$\blacksquare$

\subsubsection{$M_{i+1}$ definition}

\paragraph{Definitions:}
\begin{enumerate}
\item Define the model $M_{i+1}$ to be the minimal model within $M_{i}$
that holds axioms 1-3 above. Where the minimum is taken using the
numbering of elements of $M_{1}$ in $V$.
\item As every variable in $var_{i+1}$ has an interpretation in $M_{i+1}$
and as $var_{i}$ were elements of $M_{i}$ and as $var_{i}\subset var_{i+1}$
define the function 
\[
I_{i}:var_{i}\rightarrow M_{i+1}
\]
 to be the mapping between $var_{i}$ as elements of $M_{i}$ and
the corresponding elements in $M_{i+1}$.
\item As $M_{i+1}$ is a set model of $M_{i}$ by \ref{subsec:basic-properties-of-extension}
we know that $f_{i}:\omega\left(M_{i}\right)\rightarrow\omega\left(M_{i+1}\right)$
exist s.t $M_{i}\underset{I}{\overset{f_{i}}{\longrightarrow}}M_{i+1}$
. Define $f_{i}$ to be that function (i.e $f_{i}$ maps the omega
of $M_{i}$ to the omega of $M_{i+1}$).
\end{enumerate}

\subsubsection{Knowledge in $M_{i}$\label{subsec:Knowledge in M_=00007Bi=00007D}}

Let $\ensuremath{\phi\left(n,k\right)}$ be a two variable formula
(in $V$) which defines a function in $ZFC$ and holds the propriety
condition \ref{propriety-condition-definition} in $V\left[G_{0}\right]$.

\paragraph{Lemma:}

Within $M_{i}$ denote the axioms 1+2+3 (as interpreted in $M_{i}$)
by the set $T$ and let $C$ be the size of the universal TM. Then
\[
10k>\mathcal{L}\left(T\right)+C
\]
 where $k$ is 
\[
k=\max_{k'\in\omega\left(M_{i}\right)}\left\{ M_{i}\vDash k'-con\left(ZFC\right)\right\} 
\]

\paragraph{Proof:}

Recall the composition of $T$: 
\begin{itemize}
\item The axioms of $ZFC$ can be coded using a number $C_{1}\in\omega\left(V\right)$
(i.e $C_{1}$ is a coding of a TM that identifies the axioms of $ZFC$).
\item The axioms of $Formulas_{i+2}$ being a finite set in $V$ can be
coded using a number $C_{2}\in\omega\left(V\right)$.
\item for a given $k'\in\omega$ the axioms of $k'-con\left(ZFC\cup Formulas_{i+1}\right)$
can be coded using a number $C_{3}\in\omega\left(V\right)$. This
is a coding of a TM which takes two inputs $n',k'$ and returns true
if $n$ is an axioms of $k'-con\left(ZFC\cup Formulas_{i+1}\right)$ 
\item The number $k$ is a natural number in $\omega\left(M_{i}\right)$.
\end{itemize}
In order to identify axioms of $T$ one needs the coding of $C_{1},C_{2},C_{3}$
above and $k$ defined in axioms 2. As all $C_{1},C_{2},C_{3}$ were
standard numbers w.r.t $V$ and $k$ was non standard numbers w.r.t
$V$ we get that 
\[
C_{1}<k\,\wedge\,C_{2}<k\,\wedge\,C_{3}<k
\]
 and thus the coding of $T$ is smaller than 
\[
C_{1}+C_{2}+C_{3}+k\leq k+k+k+k=4k
\]
As $\mathcal{L}\left(T\right)$ is smaller than $2\left|T\right|$
(recall \ref{subsec:Chaitin's-incompleteness-theorem} first remark)
we know that 
\[
\mathcal{L}\left(T\right)<8k
\]
 and as $C$ the size of the universal TM is also a standard number
w.r.t $V$ we get that $C<k$ and as such 
\[
\mathcal{L}\left(T\right)+C<9k<10k
\]
 which is the assertion of the lemma. $\blacksquare$

\paragraph{Theorem:}

$M_{i}$ hold this axiom: For every $j\in\omega\left(V\right)$ it
holds 
\[
j-con\left(ZFC\cup Formulas_{i+1}\cup\varphi\right)
\]
where $\varphi$ denote the axiom ``$\exists y,c\in\omega$ s.t $y\leq10k$
and $y$ codes a TM which computes for every $n<\left(G_{0}\left(M_{i}\right)\right)_{3}$
the value of $k'\in\omega$ s.t $\phi\left(n,k'\right)$ in $\left\lceil \log\left(n+1\right)+\log\left(k'+1\right)+c\right\rceil $
steps'' and $k$ is 
\[
k=\max_{k'\in\omega\left(M_{i}\right)}\left\{ M_{i}\vDash k'-con\left(ZFC\right)\right\} 
\]
.

\paragraph{Remark:}

Please note that the statement $\varphi$ in the lemma isn't a formula
in $V$ (as it has the value of $k$ and of $\left(G_{0}\left(M_{i+1}\right)\right)_{3}$
in it, and both aren't a standard number w.r.t $V$) and as such it
and its negation may not appear in $Formulas_{i'}$ for any $i'\in\omega\left(V\right)$.
However, we still may ask the question of ``$\varphi$ is $j$ consistent
with the previous statements or not?'' which the above theorem answers.

\paragraph{Proof:}

First by previous lemma we know that axioms 1-3 are consistent in
$M_{i}$ as $M_{i}$ held 
\[
j-con\left(ZFC\cup Formulas_{i+1}\right)
\]
for every $j\in\omega\left(V\right)$. Denote $T$ to be the set of
axioms 1-3 in $M_{i}$. Let $t=\left(G_{0}\left(M_{i}\right)\right)_{3}\in\omega\left(M_{i}\right)$,
by the binary tree construction \ref{subsec:binary-tree-construction}
we know that exist $x\in\omega\left(M_{i}\right)$ codes a TM which
computes the mapping $n\rightarrow k''$ s.t $\phi\left(n,k''\right)$
holds for every $n<t$ in $\left\lceil \log\left(n+1\right)+\log\left(k''+1\right)+1\right\rceil $
computational steps. Let $k$ be the number above. 

Let $\overset{\,o}{Z}=10k$ and let $C$ be the size of the universal
TM, then we know by previous lemma that $\overset{\,o}{Z}>\mathcal{L}\left(T\right)+C$.
By reduction of machine number \ref{subsec:reduction-of-machine},
we know that a model $N$ in $M'$ of axioms 1-3 exists that holds
$\varphi$ of the axiom . Thus axioms 1-3 and $\varphi$ have a model
in $M'$ and therefore are consistent. Let $j\in\omega\left(V\right)$,
as $N$ is a model of 
\[
j-con\left(ZFC\cup Formulas_{i+1}\right)
\]
 we get that exist in $N$ a sequence of models $N_{1},\ldots,N_{j}$
each is a set of its previous and all holds $ZFC\cup Formulas_{i+1}$.
As $\varphi$ of axiom is absolute (the numbers $y,c\in\omega\left(N\right)$
exist also in $N_{i'}$ for $1\leq i'\leq j$ and holds the same condition
due to the absoluteness to TM in sub-model \ref{subsec:absoluteness-of-Turing}
and the fact that $\phi$ holds the propriety condition \ref{propriety-condition-definition}).
, As such, $N$ holds for every $j\in\omega\left(V\right)$ 
\[
j-con\left(ZFC\cup Formulas_{i+1}\cup\varphi\right)
\]
 Thus $M'$ holds 
\[
j-con\left(ZFC\cup Formulas_{i+1}\cup\varphi\right)
\]
 with the models $N_{1},\ldots,N_{j}$. As $j\in\omega\left(V\right)$
was chosen arbitrarily the statement holds for every $j\in\omega\left(V\right)$
and therefore $M_{i}$ holds the axiom.$\blacksquare$

\subsubsection{Lattice definition}

As $M_{i}$ was constructed for every $i\in\omega\left(V\right)$,
the set $\left\{ M_{i}\right\} _{i\in\omega\left(V\right)}$ exists
in $V$. 

Please note that indeed every $M_{i}$ separately is a set of $M_{1}$
(and moreover is a set of every $M_{t}$ for $t<i$ ) but as $M_{1}$
doesn't have access to $\omega\left(V\right)$ the collection $\left\{ M_{i}\right\} _{i\in\omega\left(V\right)}$
can't be defined in $M_{1}$.

\paragraph{Definition:}

The lattice is the following collection defined in $V\left[G_{0}\right]$:
\begin{itemize}
\item The set of models $\left\{ M_{i}\right\} _{i\in\omega\left(V\right)}$ 
\item The set of variables $\left\{ var_{i}\right\} _{i\in\omega\left(V\right)}$
within each model.
\item The set of formulas $\left\{ Formulas_{i}\right\} _{i\in\omega\left(V\right)}$
of each model.
\item The mappings $\left\{ I_{i}\right\} _{i\in\omega\left(V\right)}$
from variables of one model to the next.
\item The mappings $\left\{ f_{i}\right\} _{i\in\omega\left(V\right)}$
from $\omega$ of one model to the next.
\end{itemize}

\subsubsection{Definition of the lattice within $V$ \label{Definition of the lattice within V}}

\paragraph{Theorem:}

Let $k\in\omega\left(V\right)$, the lattice up to $k$ (i.e the set
$\left(M_{i},var_{i},Formulas_{i},I_{i},f_{i}\right)_{1\leq i\leq k}$)
can be defined within $V$.

\paragraph{Proof:}

Please note that for $i\in\omega\left(V\right)$ the construction
above of $M_{i+1},I_{i},f_{i},var_{i+1},Formulas_{i+1}$ doesn't depend
on $G_{0}$ entirely. It only depends on the values of $G_{0}\left(M_{i}\right)$.
In other words, in order to define $M_{i+1}$ we don't need to know
all of $G_{0}$ but we only need to know a partial function on $P_{0}$
which contains $\left\{ M_{k}\,|\,k\leq i\right\} $ in its domain.
As the restriction of $G_{0}$ to the set $\left\{ M_{k}\,|\,k\leq i\right\} $
gives us a partial function in $P_{0}$ and as such partial function
are assumed to be in $V$ we get that the set $\left(M_{i},var_{i},Formulas_{i},I_{i},f_{i}\right)_{1\leq i\leq k}$
can be defined within $V$. $\blacksquare$

\subsubsection{independence of the lattice from $\left(G_{0}\right)_{3}$ \label{independence of the lattice from (G_=00007B0=00007D)_=00007B3=00007D}}

\paragraph{Observation:}

Recall that the demands on $M_{i+1}$ were axioms 1-3 and recall that
axioms 1-3 didn't use the values of $\left(G_{0}\right)_{3}$. Therefore
the lattice is constructed independently of the values $\left(G_{0}\right)_{3}$.
Or, in other words, if the values in $\left(G_{0}\right)_{3}$ were
to be change the lattice under the above construction would remain
the same. 

\section{Limit of the lattice\label{sec:Limit-of-the-lattice}}

In this part we will use the construction of the lattice defined above
to create a ``limit'' model $M_{Limit}$ of models models $\left\{ M_{i}\right\} _{i\in\omega\left(V\right)}$
as a type of straight limit over the sets $\left\{ var_{i}\right\} _{i\in\omega\left(V\right)}$
and the mappings $\left\{ I_{i}\right\} _{i\in\omega\left(V\right)}$
in such a way that, the limit model will hold $ZFC$ and will have
good knowledge properties (we will define them in \ref{subsec:Knowledge-in-M_limit})

\subsection{Limit model definition\label{subsec:Limit-model-definition}}

\paragraph{Definition:}

Denote the set 
\[
Lace_{i}=\left\{ \left(a_{k}\right)_{k\in\omega\left(V\right)\wedge k\geq i}\,|\,\begin{matrix}a_{i}\in var_{i}\wedge\\
\forall k\geq i\,\,a_{k+1}=I_{k}\left(a_{k}\right)
\end{matrix}\right\} 
\]
 In other words, as every element in $var_{i}$ is mapped into $var_{i+1}$
we denote their trajectories by $Lace_{i}$.

\paragraph{Definition:}

Denote 
\[
Lace=\cup_{i\in\omega\left(V\right)}Lace_{i}
\]

\paragraph{Lemma:}

If $a\in Lace_{i_{1}}$ and $b\in Lace_{i_{2}}$ let $k\geq\max\left\{ i_{1},i_{2}\right\} $
s.t $a_{k}=b_{k}$. Then for all $k'\geq k$ it holds $a_{k'}=b_{k'}$.

\paragraph{Proof:}

By induction over $k'\geq k$: The base case is $a_{k}=b_{k}$ and
is assumed true. Assume $a_{k'}=b_{k'}$ for $k'\geq k$ then $a_{k'+1}=I_{k'}\left(a_{k'}\right)=I_{k'}\left(b_{k'}\right)=b_{k'+1}$
$\blacksquare$

\paragraph{Definition:}

Define the following equivalence relation on elements of $Lace$:
for $a\in Lace_{i_{1}}$ and $b\in Lace_{i_{2}}$ we denote $a\sim b$
if exist $k\geq\max\left\{ i_{1},i_{2}\right\} $ s.t $a_{k}=b_{k}$.

\subsubsection{Limit model set definition on }

\paragraph{Definition:}

Define 
\[
M_{Limit}=Lace/\sim
\]
 the equivalence class of $Lace$ under $\sim$ above.

\subsubsection{Limit model $\in$ definition \label{Limit-model-=00005Cin-definition}}

\paragraph{Definition:}

For $a\in Lace_{i_{1}}$ and $b\in Lace_{i_{2}}$ we denote $a\in_{Lace}b$
if exist $k\geq\max\left\{ i_{1},i_{2}\right\} $ s.t $\forall k'\geq k$
it holds that $a_{k'}\in_{M_{k'}}b_{k'}$.

\paragraph{Definition:}

For $a,b\in M_{Limit}$ define $a\in_{Limit}b$ if for all $a',b'$
being representatives of the equivalent classes $a,b$ respectively
it holds $a'\in_{Lace}b'$.

\paragraph{Remark:}

In the following lemma we will prove that $a\in_{Limit}b$ is independent
of the choice of representatives, so this definition could have been
define as exist $a',b'$ representatives s.t $a'\in_{Lace}b'$. 

\paragraph{Lemma:}

Let $a,b\in M_{Limit}$ then $a\in_{Limit}b$ if and only if exist
$a'',b''$ representatives of the equivalent classes $a,b$ respectively
that holds $a''\in_{Lace}b''$.

\paragraph{Proof:}

First assume that exist $a'',b''$ representatives s.t $a''\in_{Lace}b''$
and let $a',b'$ be another pair of representatives of $a,b$. Then
(as they represent the same equivalent class) $a'\sim a''$ and $b'\sim b''$.
By the definition of $\sim$ exist $k_{1},k_{2}\in\omega\left(V\right)$
s.t $\forall k>k_{1}\,\,a'_{k}=a''_{k}$ and $\forall k>k_{2}\,\,b'_{k}=b''_{k}$.
As $a''\in_{Lace}b''$ by definition exist $k_{3}\in\omega\left(V\right)$
s.t $\forall k>k_{3}\,\,a''_{k}\in_{M_{k}}b''_{k}$ denote $k'=\max\left\{ k_{1},k_{2},k_{3}\right\} $
then for $k>k'$ it holds that 
\[
a'_{k}=a''_{k}\in_{M_{k}}b''_{k}=b'_{k}\Rightarrow a'{}_{k}\in_{M_{k}}b'{}_{k}
\]
 and therefore $a'\in_{Lace}b'$ and hence as $a',b'$ were any pair
of representatives of $a,b$ we get $a\in_{Limit}b$ .

The other direction is trivial: If $a\in_{Limit}b$ then by definition
exist $a',b'$ representatives of $a,b$ and as $a''\in_{Lace}b''$
for any two representatives $a'',b''$ we get $a'\in_{Lace}b'$ for
the case of $a',b'$ as well. $\blacksquare$

\subsection{Stabilization theorems}

\subsubsection{Appearance of formulas in $Formulas_{k}$ \label{Appearance-of-formulas in Formulas_=00007Bk=00007D }}

\subparagraph{Theorem:}

Let $\psi\left(x_{1},\ldots x_{n}\right)$ be a formula of the language
of set theory (in $V$) with $n\in\omega\left(V\right)$ free variables
and let $a_{1},\ldots,a_{n}\in Lace$ be elements of $Lace$ then
exist $k\in\omega\left(V\right)$ s.t either the formula $\psi\left(a_{1,k},\ldots a_{n.k}\right)$
or it negation is in $Formulas_{k}$. Where $a_{i,k}$ is the element
$a_{i}$ at $var_{k}$.

\subparagraph{Note:}

Please note that in the above theorem we didn't argue that the truth
value of $\psi$ is stabilized from a certain point on, this will
be proven later. This theorem shows that the formula $\psi$ itself
is to be found within $Formulas_{k}$.

\subparagraph{Proof:}

Let $k\in\omega\left(V\right)$ be large enough s.t all $a_{1,i_{1}},\ldots,a_{n,i_{1}}$
are to be found in $var_{k}$. Recall that by \ref{Definition of the lattice within V}
$k\in\omega\left(V\right)$ we can define $\left(M_{i},var_{i},Formulas_{i},I_{i},f_{i}\right)_{1\leq i\leq k}$
within $V$. Define the following set $E$ in $V$: $E$ is the set
of partial function $f$ in $P_{0}$ that agreed with $G_{0}$ on
$\left(M_{i}\right)_{1\leq i\leq k}$ and such that if we continue
the construction of the lattice according to $f$ we will get the
formula $\psi$ (or it negation) in some later construction of $Formulas_{k'}$
according to $f$. As every partial function that extends $G_{0}\restriction_{\left(M_{i}\right)_{1\leq i\leq k}}$
can be extended into a function in $E$ we get that $E$ is a set
in $V$ that is dense over $G_{0}\restriction_{\left(M_{i}\right)_{1\leq i\leq k}}$
and thus by \ref{subsec:Forcing} theorem 3 we get that $E\cap G_{0}\neq\emptyset$.
That is, exist $k'\in\omega\left(V\right)$ s.t $\psi$ (or its negation)
appear in $Formulas_{k'}$ (and now the construction is done according
to $G_{0}$, the usual way). $\blacksquare$

\subsubsection{Stabilization in $Formulas_{k}$ \label{Stabilization in Formulas_=00007Bk=00007D}}

\subparagraph{Lemma:}

Let $\psi\left(x_{1},\ldots x_{n}\right)$ be a formula of the language
of set theory (in $V$) with $n\in\omega\left(V\right)$ free variables
and let $a_{1},\ldots,a_{n}\in Lace$ be elements of $Lace$ then
exist $k\in\omega\left(V\right)$ s.t exactly one of the following
holds:
\begin{enumerate}
\item for all $k'\in\omega\left(V\right)$ s.t $k'>k$ the formula $\psi\left(a_{1,k'},\ldots a_{n.k'}\right)$
is in $Formulas_{k'}$.
\item for all $k'\in\omega\left(V\right)$ s.t $k'>k$ the formula $\neg\psi\left(a_{1,k'},\ldots a_{n.k'}\right)$
is in $Formulas_{k'}$.
\end{enumerate}

\paragraph{Terminology remark:}

If for all $k'\in\omega\left(V\right)$ s.t $k'>k$ the formula $\psi\left(a_{1,k},\ldots a_{n.k}\right)$
is in $Formulas_{k'}$ we say that the formula $\psi$ had \textbf{stabilized}.
This is different than the assertion in \ref{Appearance-of-formulas in Formulas_=00007Bk=00007D },
as in that section we only argued that the formula $\psi$ or $\neg\psi$
\textbf{appears} in $Formulas_{k'}$, but perhaps it may be the case
that for even $k'$s $\psi$ is in $Formulas_{k'}$ and for odd $k'$s
$\neg\psi$ is in $Formulas_{k'}$? Stabilization is the assertion
that such phenomena doesn't happen as either $\psi$ appears in $Formulas_{k'}$
from one point on or $\neg\psi$ appears in $Formulas_{k'}$ from
one point on.

\paragraph{Proof:}

By previous theorem \ref{Appearance-of-formulas in Formulas_=00007Bk=00007D }
we know that exist $k_{0}\in\omega\left(V\right)$ s.t either the
formula $\psi\left(a_{1,k_{0}},\ldots a_{n.k_{0}}\right)$ or it negation
is in $Formulas_{k_{0}}$. Recall that $Formulas_{k}$ is a function
(i.e an ordered set of formulas) the proof is by induction over the
place number of $\psi$ in the list. Recall the definition of $Formulas_{k_{0}+1}$:
\[
Formulas_{k_{0}+1}=Reduce_{M_{k_{0}}}\left(Reduce_{M_{k_{0}}}\left(\emptyset,Formulas_{k_{0}}\right),\left(G_{3}\left(M_{k_{0}}\right)\right)_{2}\right)
\]
and as such every $\psi$ that appear in $Formulas_{k_{0}}$ must
appear in $Formulas_{k_{0}+1}$ as the formula itself or its negation.
By induction we receive that for all $k'>k_{0}$ the formula $\psi$
or it negation appear in $Formulas_{k'}$. And moreover, the if the
formula or its negation held place $i$ in $Formulas_{k_{0}}$, the
same formula or its negation will have place $i$ in $Formulas_{k'}$.
Now we will prove the lemma by induction over the place number of
$\psi$ in the list.

\subparagraph{Base case: }

The formula $\psi$ appear as first formula in the list:

If for all $k'>k_{0}$ the formula $\psi$ (and not $\neg\psi$) appears
in $Formulas_{k'}$ then the lemma holds for $k=k_{0}$ with condition
(1). Otherwise exist $k_{1}>k_{0}$ s.t $\neg\psi$ appear in $Formulas_{k_{1}}$.
By the definition of $Reduce_{M_{k_{1}}}$ it means that the formula
$\psi$ can be proven false from a finite (by $M_{k_{1}}$) subset
of the axioms of 
\[
\left\{ s-con\left(ZFC\right)\,|\,s\in\omega\left(V\right)\right\} 
\]
For $k'>k_{1}$, as the set 
\[
\left\{ s-con\left(ZFC\right)\,|\,s\in\omega\left(V\right)\right\} 
\]
doesn't change and as every proof in $M_{k_{1}}$ is also a proof
in $M_{k'}$ it must be the case that $\neg\psi$ appear in $Formulas_{k'}$
as well. In this case, the lemma holds for $k=k_{1}$ with condition
(2).

\subparagraph{Step case:}

The formula $\psi$ appear as $t+1$ formula in the list for $t\in\omega\left(V\right)$. 

In this case let by the induction's assumption we know that for every
formula up to $t$ exist $k'_{t}$ s.t for all $k'>k_{t}$ the formula
stabilizes in $Formulas_{k'}$. Denote $k_{0}=\max\left\{ k'_{0},\ldots,k'_{t}\right\} $.
We know that all formulas in $Formulas_{k_{0}}$ up to $\psi$ had
stabilized. If for all $k'>k_{0}$ the formula $\psi$ (and not $\neg\psi$)
appears in $Formulas_{k'}$ then the lemma holds for $k=k_{0}$ with
condition (1). Otherwise exist $k_{1}>k_{0}$ s.t $\neg\psi$ appear
in $Formulas_{k_{1}}$. By the definition of $Reduce_{M_{k_{1}}}$
it means that the formula $\psi$ can be proven false from a finite
(by $M_{k_{1}}$) subset of the axioms of 
\[
\left\{ s-con\left(ZFC\cup A\right)\,|\,s\in\omega\left(V\right)\right\} 
\]
 where $A$ is the set of formulas that appear before $\psi$ and
$t\in\omega\left(V\right)$. For $k'>k_{1}$, as the set 
\[
\left\{ s-con\left(ZFC\cup A\right)\,|\,s\in\omega\left(V\right)\right\} 
\]
doesn't change, as all formulas in $A$ stabilized up to $z_{0}$
and as every proof in $M_{k_{1}}$ is also a proof in $M_{k'}$ it
must be the case that $\neg\psi$ appear in $Formulas_{k'}$ as well.
In this case, the lemma holds for $k=k_{1}$ with condition (2). $\blacksquare$

\paragraph{Remark:}

As it may be unclear, for a set $B$ of formulas that exist in $M_{i}$
if $B$ is inconsistent, then a contradiction can be proven from a
finite set of statements in $B$. The same proof is valid in any set
model of $M_{i}$ as proofs are verifiable using a TM and due to absoluteness
of Turing machines in sub-models \ref{subsec:absoluteness-of-Turing}.

\subsubsection{Stabilization of formulas in $Lace$}

\paragraph{Theorem:}

Let $n\in\omega\left(V\right)$ and $a_{1},\ldots,a_{n}\in Lace$
s.t $a_{k}\in Lace_{i_{k}}$. Let $\psi$ be a formula (in $V$) of
the language of set theory with $n$ free variables. then exist $k\in\omega\left(V\right)$,
$k\geq\max\left\{ i_{1},i_{2},\ldots,i_{n}\right\} $ s.t exactly
one of the following holds:
\begin{enumerate}
\item for all $k'\in\omega\left(V\right)$ s.t $k'>k$ $M_{k'}\models\psi\left(a_{1,k},\ldots,a_{n,k}\right)$
.Where $a_{i,k}$ is the element $a_{i}$ at $var_{k}$.
\item for all $k'\in\omega\left(V\right)$ s.t $k'>k$ $M_{k'}\models\neg\psi\left(a_{1,k},\ldots,a_{n,k}\right)$
.Where $a_{i,k}$ is the element $a_{i}$ at $var_{k}$.
\end{enumerate}

\paragraph{Proof:}

By \ref{Stabilization in Formulas_=00007Bk=00007D} we know that exist
$k_{0}\in\omega\left(V\right)$ s.t $\phi$ stabilizes from the $k_{0}$
place onward on $Formulas_{k}$. Denote $k=k_{0}+1$. 

If for $k'\in\omega\left(V\right)$ s.t $k'>k_{0}$ the formula $\psi\left(a_{1,k'},\ldots a_{n.k'}\right)$
is in $Formulas_{k'}$. Then, as $k'\geq k_{0}+1$ we know that formula
$\psi\left(a_{1,k'-1},\ldots a_{n.k'-1}\right)$ is in $Formulas_{k'-1}$
and as such all formulas in $Formulas_{k'-1}$ are axioms of $M_{k'}$
if follows that $M_{k'}\models\psi\left(a_{1,k'},\ldots,a_{n,k'}\right)$
and the theorem holds with condition (1). 

Otherwise, If for $k'\in\omega\left(V\right)$ s.t $k'>k_{0}$ the
formula $\neg\psi\left(a_{1,k'},\ldots a_{n.k'}\right)$ is in $Formulas_{k'}$
Then, as $k'\geq k_{0}+1$ we know that formula $\neg\psi\left(a_{1,k'-1},\ldots a_{n.k'-1}\right)$
is in $Formulas_{k'-1}$ and as such all formulas in $Formulas_{k'-1}$
are axioms of $M_{k'}$ if follows that $M_{k'}\models\neg\psi\left(a_{1,k'},\ldots,a_{n,k'}\right)$
and the theorem holds with condition (2). $\blacksquare$

\subsubsection{Witnesses of formulas in $Lace$\label{Witnesses-of-formulas-in-Lace}}

\paragraph{Theorem:}

Let $n\in\omega\left(V\right)$ and $a_{1},\ldots,a_{n}\in Lace$
s.t $a_{k}\in Lace_{i_{k}}$. Let $\psi$ be a formula (in $V$) of
the language of set theory with $n+1$ free variables. Assume that
exist $k\in\omega\left(V\right)$, $k\geq\max\left\{ i_{1},i_{2},\ldots,i_{n}\right\} $
s.t for all $k'\in\omega\left(V\right)$ s.t $k'>k$ $M_{k'}\models\exists y\,\psi\left(y,a_{1,k'},\ldots,a_{n,k'}\right)$.
Then, exist $z\in Lace$ and $k''>k$ s.t for all $k'\in\omega\left(V\right)$
s.t $k'>k''$ $M_{k'}\models\psi\left(z_{k'},a_{1,k'},\ldots,a_{n,k'}\right)$.

\paragraph{Proof:}

Let $k_{0}\in\omega\left(V\right)$ be a number s.t $\exists y\,\psi\left(y,a_{1,k'},\ldots,a_{n,k'}\right)$
stabilizes and all previous formulas up to it had stabilized as well.
Let the set $E$ be the set of partial functions $g=\left(g_{1},g_{2}\right)$
extending $G_{0}$ on $\left\{ M_{i}\right\} _{i=2}^{k_{0}}$ in $P_{0}$
s.t for every finite sequence of models $\left\{ M'_{i}\right\} _{i\in\left[a\right]}$
($a\in\omega\left(V\right))$ in $Dom\left(g\right)$ with $M'_{0}=M_{k_{0}}$
if the construction of the lattice was made using 
\[
M_{1},\ldots,M_{k_{0}}=M'_{0},M'_{1},\ldots,M'_{a}
\]
 then exists $i\leq a$ s.t if $M'_{i}\models\exists y\,\psi\left(y,a'_{1,i},\ldots,a'_{n,i}\right)$
with $a'_{1,k'},\ldots,a'_{n,k'}$ the appropriate variables in $var$
of $M'_{i}$ then $y'$ is chosen to be a variable in $g_{1}$ that
isn't in $Var_{i-1}$ and the formula 
\[
\psi\left(y',a'_{1,k'},\ldots,a'_{n,k'}\right)
\]
 was chosen as a formula in $g_{2}$.

As $Var_{i-1}$ is a finite (by $V$) set there is $y'$ that wasn't
chosen. partial function can be extended to include $y'$ and the
formula $\psi$, the set $E$ is dense above the partial function
of $G_{3}$ on $\left\{ M_{i}\right\} _{i=2}^{k_{0}}$ in $V$. So
it must be the case that $E\cap G_{0}\neq\emptyset$.

Therefore, exist $k_{1}>k_{0}$ s.t an element $y'$ exists in $var_{k_{1}}$
but not in $var_{k_{1}-1}$ and the formula
\[
\psi\left(y'_{k_{1}},a'_{1,k'},\ldots,a'_{n,k'}\right)
\]
 is in $Formulas_{k_{1}}$ . As consistency, $y'$ didn't appear in
$var_{k_{1}-1}$ so the above formula is equivalent to $\exists y\,\psi\left(y,a'_{1,i},\ldots,a'_{n,i}\right)$
(as $y'$ may be mapped to any element, as no other axioms are demanded
on it). We know that the formula $\exists y\,\psi\left(y,a'_{1,i},\ldots,a'_{n,i}\right)$
is consistent with the previous formulas (as it holds true from $M_{k_{1}+1}$
and the latter is a set of $M_{k_{1}}$). Let $k'>k_{1}$, as $M_{k'+1}\models\exists y\,\psi\left(y,a_{1,k'},\ldots,a_{n,k'}\right)$
we get that $\exists y\,\psi\left(y,a'_{1,i},\ldots,a'_{n,i}\right)$
is consistent with the previous formulas and as $y'$ appear first
in $\psi\left(y'_{k'},a_{1,k'},\ldots,a_{n,k'}\right)$ we know that
it is equivalent to $\exists y\,\psi\left(y,a'_{1,i},\ldots,a'_{n,i}\right)$
hence $M_{k'}\models\psi\left(y'_{k'},a_{1,k'},\ldots,a_{n,k'}\right)$.
Let $z=y'_{k_{1}},y'_{k_{1}+1},y'_{k_{1}+2},\ldots$ be that element
in $Lace$ and $k''=k_{1}$. We get that for all $k'\in\omega\left(V\right)$
s.t $k'>k''$ $M_{k'}\models\psi\left(z_{k'},a_{1,k'},\ldots,a_{n,k'}\right)$.
$\blacksquare$

\subsubsection{Truth value in $M_{Limit}$ \label{4.2.5 Truth value in M_=00007BLimit=00007D }}

\paragraph{Theorem:}

Let $n\in\omega\left(V\right)$ and $a_{1},\ldots,a_{n}\in M_{Limit}$.
Let $\psi$ be a formula (in $V$) of the language of set theory with
$n$ free variables. Then $M_{Limit}\models\psi\left(a_{1},\ldots,a_{n}\right)$
if and only if exist $a_{1},\ldots,a_{n}\in Lace$ s.t $a_{k}\in Lace_{i_{k}}$
which are representatives of $a_{1},\ldots,a_{n}\in M_{Limit}$ respectively.
And exist $k\in\omega\left(V\right)$, $k\geq\max\left\{ i_{1},i_{2},\ldots,i_{n}\right\} $
s.t for all $k'\in\omega\left(V\right)$ s.t $k'>k$ $M_{k'}\models\psi\left(a_{1,k'},\ldots,a_{n,k'}\right)$.
Where $a_{i,k}$ is the element $a_{i}$ at $var_{k}$.

\paragraph{Remark:}

As the truth value doesn't depend on the representative we could have
said ``... if and only if for all $a_{1},\ldots,a_{n}\in Lace$ s.t
$a_{k}\in Lace_{i_{k}}$ which are representatives of $a_{1},\ldots,a_{n}\in M_{Limit}$
respectively...''.

\paragraph{Proof:}

The proof is by induction over the structure of the formula $\psi$
where the formula is in $V$ and the induction is done in $V$ as
well:

Base cases:
\begin{itemize}
\item For the basic case $\psi\left(a',b'\right)="a\in b"$ if for all $k'\in\omega\left(V\right)$
s.t $k'>k$ $M_{k'}\models\left(a_{k'}\in_{M_{k'}}b_{k'}\right)$
then by definition of $\in_{Limit}$ (see \ref{Limit-model-=00005Cin-definition})
it holds $M_{Limit}\models a\in_{Limit}b$. On the other hand, if
$M_{Limit}\models a\in_{Limit}b$ then by definition exist $k\geq\max\left\{ i_{1},i_{2}\right\} $
s.t $\forall k'\geq k$ it holds that $a_{k'}\in_{M_{k'}}b_{k'}$.
\item For the basic case of $\psi\left(a',b'\right)="a=b"$ if exist $k'\in\omega\left(V\right)$
s.t $M_{k'}\models\left(a_{k'}=b_{k'}\right)$ the mapping $I_{i}$
will map both $a_{k'},b_{k'}$ to the same element in $M_{k'+1}$
and therefore both laces $a',b'$ are equivalent under $\sim$ and
as such represent the same element in $M_{Limit}$ (i.e the are in
the same equivalence class) and as such $M_{Limit}\models\left(a=b\right)$.
On the other hand, if $M_{Limit}\models\left(a=b\right)$ then by
definition of equality $a,b$ represent the same element in $M_{Limit}$
choose $a',b'\in Lace$ which represent that (same) element and as
such $M_{k'}\models\left(a_{k'}=b_{k'}\right)$ for $k'$ large enough
and on.
\end{itemize}
Composite case:
\begin{itemize}
\item For the case of $\psi\left(a_{1},\ldots,a_{n}\right)=\psi_{1}\left(a_{1},\ldots,a_{n}\right)\wedge\psi_{2}\left(a_{1},\ldots,a_{n}\right)$,
as exist $k\in\omega\left(V\right)$, $k\geq\max\left\{ i_{1},i_{2},\ldots,i_{n}\right\} $
s.t for all $k'\in\omega\left(V\right)$ s.t $k'>k$ $M_{k'}\models\psi\left(a'_{1,k'},\ldots,a'_{n,k'}\right)$
it follows that $M_{k'}\models\psi_{1}\left(a'_{1,k'},\ldots,a'_{n,k'}\right)$
and $M_{k'}\models\psi_{2}\left(a'_{1,k'},\ldots,a'_{n,k'}\right)$.
As $\psi_{1},\psi_{2}$ were formulas of lower depth we get from the
induction's assumption that $M_{Limit}\models\psi_{1}\left(a_{1},\ldots,a_{n}\right)$
and $M_{Limit}\models\psi_{2}\left(a_{1},\ldots,a_{n}\right)$ and
therefore $M_{Limit}\models\psi\left(a_{1},\ldots,a_{n}\right)$.
On the other hand, if $M_{Limit}\models\psi\left(a_{1},\ldots,a_{n}\right)$
then we know that $M_{Limit}\models\psi_{1}\left(a_{1},\ldots,a_{n}\right)$
and $M_{Limit}\models\psi_{2}\left(a_{1},\ldots,a_{n}\right)$ and
therefore exists $k_{1},k_{2}$ s.t for all $k'\in\omega\left(V\right)$
s.t $k'>\max\left\{ k_{1},k_{2}\right\} $ it holds $M_{k'}\models\psi_{1}\left(a'_{1,k'},\ldots,a'_{n,k'}\right)$
and $M_{k'}\models\psi_{2}\left(a'_{1,k'},\ldots,a'_{n,k'}\right)$
it follows that $M_{k'}\models\psi\left(a'_{1,k'},\ldots,a'_{n,k'}\right)$
as needed. 
\item The case $\psi\left(a_{1},\ldots,a_{n}\right)=\psi_{1}\left(a_{1},\ldots,a_{n}\right)\vee\psi_{2}\left(a_{1},\ldots,a_{n}\right)$
is analogues to the previous case.
\item For the case $\psi\left(a_{1},\ldots,a_{n}\right)=\neg\psi_{1}\left(a_{1},\ldots,a_{n}\right)$,
as $\psi_{1}$ formula of lower depth we get from the induction's
assumption that exist $k\in\omega\left(V\right)$, $k\geq\max\left\{ i_{1},i_{2},\ldots,i_{n}\right\} $
s.t for all $k'\in\omega\left(V\right)$ s.t $k'>k$ $M_{k'}\models\psi_{1}\left(a'_{1,k'},\ldots,a'_{n,k'}\right)$
if and only if $M_{Limit}\models\psi_{1}\left(a_{1},\ldots,a_{n}\right)$.
As such, the same assertion is true of the negation of the formula,
we get that exist $k\in\omega\left(V\right)$, $k\geq\max\left\{ i_{1},i_{2},\ldots,i_{n}\right\} $
s.t for all $k'\in\omega\left(V\right)$ s.t $k'>k$ $M_{k'}\models\neg\psi_{1}\left(a'_{1,k'},\ldots,a'_{n,k'}\right)$
if and only if $M_{Limit}\models\neg\psi_{1}\left(a_{1},\ldots,a_{n}\right)$. 
\end{itemize}
quantifier case:
\begin{itemize}
\item For the case $\psi\left(a_{1},\ldots,a_{n}\right)=\exists y\,\psi_{1}\left(y,a_{1},\ldots,a_{n}\right)$:
If exist $k\in\omega\left(V\right)$, $k\geq\max\left\{ i_{1},i_{2},\ldots,i_{n}\right\} $
s.t for all $k'\in\omega\left(V\right)$ s.t $k'>k$ $M_{k'}\models\exists y\,\psi_{1}\left(y,a_{1,k'},\ldots,a_{n,k'}\right)$
it follows by the existence of witnesses in $Lace$ (see \ref{Witnesses-of-formulas-in-Lace})
that exist $z'\in Lace$ and $k''>k$ s.t for all $k'\in\omega\left(V\right)$
s.t $k'>k''$ $M_{k'}\models\psi_{1}\left(z'_{k'},a_{1,k'},\ldots,a_{n,k'}\right)$.
As $\psi_{1}$ is a formula of less depth we get that $M_{Limit}\models\psi_{1}\left(z,a_{1},\ldots,a_{n}\right)$
for $z\in M_{Limit}$ being the equivalence class of $z'$. As exist
$z\in M_{Limit}$ s.t $M_{Limit}\models\psi_{1}\left(z,a_{1},\ldots,a_{n}\right)$
it holds $M_{Limit}\models\exists y\,\psi_{1}\left(y,a_{1},\ldots,a_{n}\right)$.
On the other hand, if $M_{Limit}\models\exists y\,\psi_{1}\left(y,a_{1},\ldots,a_{n}\right)$
then exist $z\in M_{Limit}$ s.t $M_{Limit}\models\psi_{1}\left(z,a_{1},\ldots,a_{n}\right)$
then as $\psi_{1}$ is a formula of less depth we know that exist
$z'\in Lace_{i_{n+1}}$ and exist $k\in\omega\left(V\right)$, $k\geq\max\left\{ i_{1},i_{2},\ldots,i_{n},i_{n+1}\right\} $
s.t for all $k'\in\omega\left(V\right)$ s.t $k'>k$ $M_{k'}\models\psi_{1}\left(z'_{k'},a_{1,k'},\ldots,a_{n,k'}\right)$
as $z'_{k}\in M_{k'}$ exist s.t $\psi_{1}\left(z'_{k'},a_{1,k'},\ldots,a_{n,k'}\right)$
we get that $M_{k'}\models\exists y\,\psi_{1}\left(y,a_{1,k'},\ldots,a_{n,k'}\right)$.
\item For the case $\psi\left(a_{1},\ldots,a_{n}\right)=\forall y\,\psi_{1}\left(y,a_{1},\ldots,a_{n}\right)$
: The case is equivalent to $\psi\left(a_{1},\ldots,a_{n}\right)=\neg\exists y\,\neg\psi_{1}\left(y,a_{1},\ldots,a_{n}\right)$
which is composed of formulas of the previous cases.
\end{itemize}

\subsubsection{$ZFC$ holds in $M_{Limit}$ \label{subsec:ZFC holds in M_=00007BLimit=00007D }}

\paragraph{Theorem:}

The set $M_{Limit}$ along with $\in_{Limit}$ holds the axioms of
$ZFC$ and as such $\left(M_{Limit},\in_{Limit}\right)$ is a $ZFC$
model.

\paragraph{Proof:}

Let $\psi$ be a formula (in $V$) which is an axioms of $ZFC$ (such
a formula must be without free variables). For all $k'\in\omega\left(V\right)$
as $M_{k'}$ is a $ZFC$ model we know $M_{k'}\models\psi$. By \ref{4.2.5 Truth value in M_=00007BLimit=00007D }
we know that if $M_{k'}\models\psi$ holds from a certain $k$ and
on then $M_{Limit}\models\psi$. As such $M_{Limit}$ holds all axioms
of $ZFC$ and as such $\left(M_{Limit},\in_{Limit}\right)$ is a $ZFC$
model. $\blacksquare$

\subsection{Knowledge in $M_{Limit}$ \label{subsec:Knowledge-in-M_limit}}

\subsubsection{$M_{Limit}$ knows $\phi$ function}

Let $\ensuremath{\phi\left(n,k\right)}$ be a two variable formula
(in $V$) which defines a function in $ZFC$ and holds the propriety
condition \ref{propriety-condition-definition} in $V\left[G_{0}\right]$.
Recall that $M_{Limit}$ was a limit of models $M_{k}$. We now are
ready to show that the model $M_{Limit}$ holds that exists $\overset{\,o}{Z}\in_{Limit}\omega\left(M_{Limit}\right)$
and for all elements $n\in_{Limit}\omega\left(M_{Limit}\right)$ exists
$x<\overset{\,o}{Z}$ and $c\in_{Limit}\omega\left(M_{Limit}\right)$
s.t $x$ represents a TM which calculates the value $k\in_{Limit}\omega\left(M_{Limit}\right)$
s.t $M_{Limit}\models\phi\left(n,k\right)$ in $\left\lceil \log\left(n+1\right)+\log\left(k+1\right)+c\right\rceil $
steps and $c$ depends only on $x$ (and not on $n$). In other words
$M_{Limit}$ will hold 
\[
\begin{matrix}\exists\overset{\,o}{Z}\in\omega\,\,\forall n\in\omega\,\,\exists x,k,c\in\omega\\
\left(\begin{matrix}\phi\left(n,k\right)\,\wedge\,Run\left(x,n\right)=k\,\,\wedge\,x<\overset{\,o}{Z}\,\,\wedge\\
\left(\forall n',k'\in\omega\right)\left(\begin{matrix}\left(\phi\left(n',k'\right)\wedge Run\left(x,n'\right)=k'\right)\rightarrow\\
\left(Time\left(x,n'\right)=\left\lceil \log\left(n'+1\right)+\log\left(k'+1\right)+c\right\rceil \right)
\end{matrix}\right)
\end{matrix}\right)
\end{matrix}
\]
 where $Run\left(x,n\right)$ denote that function that return the
value return by the TM numbered $x$ on input $n$ and $Time\left(x,n\right)$
return the number of steps done in the calculation of $x$ on input
$n$.

\subsubsection{Consistency number in $M_{Limit}$\label{subsec:Consistency-number-in-M_limit}}

\paragraph{Definition:}

Define the consistency number
\[
cn=\left(cn_{i}\right)_{i\in\omega\left(V\right)}\in\underset{i\in\omega\left(V\right)}{\Pi}\omega\left(M_{i}\right)
\]
to be, for a given $i\in\omega\left(V\right)$ the number $cn_{i}\in_{M_{i}}\omega\left(M_{i}\right)$
is the only number s.t 
\[
M\vDash\left(cn_{i}\right)-con\left(ZFC\right)
\]
 but not 
\[
M\not\vDash\left(cn_{i}+1\right)-con\left(ZFC\right)
\]

\paragraph{Remark:}

In other words $cn_{i}$ is defined as
\[
cn_{i}=\max_{k'\in\omega\left(M_{i}\right)}\left\{ M_{i}\vDash k'-con\left(ZFC\right)\right\} 
\]
 as every $M_{i}$ is a set of $M_{1}$ and as $M_{1}$ doesn't hold
$\omega-con\left(ZFC\right)$ such $cn_{i}$ must exist in every $M_{i}$ 

\paragraph{Definition:}

For $k\in\omega\left(V\right)$ denote $cn\restriction_{<k}$ to be
\[
cn\restriction_{<k}=\left(cn_{i}\right)_{\begin{matrix}i\in\omega\left(V\right)\\
i>k
\end{matrix}}\in\underset{\begin{matrix}i\in\omega\left(V\right)\\
i>k
\end{matrix}}{\Pi}\omega\left(M_{i}\right)
\]

\paragraph{Lemma:}

Exist $k\in\omega\left(V\right)$ s.t $cn\restriction_{<k}\in Lace_{k}$ 

\paragraph{Proof :}

Denote by $\psi\left(x\right)$ the formula, in the language of set
theory, that states ``$x\in\omega$ and $x-con\left(ZFC\right)$
and $\neg\left(\left(x+1\right)-con\left(ZFC\right)\right)$''. Define
the following set $E$ in $V$:

$E$ is the set of partial functions $f$ in $P_{0}$ such that if
we continue the construction of the lattice according to $f$ we will
get $k'\in\omega\left(V\right)$ s.t $k'>1$ and we will get a variable
$x$ that didn't appear in $var_{k'-1}$ and $\psi\left(x\right)$
(or its negation) appear in $Formulas_{k'}$. As every partial function
can be extended into a function in $E$ we get that $E$ is a set
in $V$ that is dense and thus by \ref{subsec:Forcing} theorem 3
we get that $E\cap G_{0}\neq\emptyset$. Therefore, a $k'\in\omega\left(V\right)$
exist s.t the variable $x$ is a new variable in $Var_{k'}$ and $\psi$
appears in $Formulas_{k'}$. Thus, as $\psi\left(x\right)$ is consistent
with the construction (as for every $M_{i}$, the model has such $x$)
the formula $\psi\left(x\right)$ (and not it's negation) appears
in $Formulas_{k'}$. Let $x_{t}$ for $t>k'$ denote the trajectory
of $x$ under $I_{t}$. Thus by definition $cn_{t}=x_{t}$ for every
$t>k'$ and as such $cn\restriction_{<k}\in Lace_{k}$ for $k=k'+1$.
$\blacksquare$

\paragraph{Corollaries:}
\begin{enumerate}
\item As $cn$ is in $Lace$ from a certain $k$ and on, there is an element
$cn\in M_{Limit}$ that is represented by $cn\restriction_{<k}\in Lace_{k}$
for the appropriate $k\in\omega\left(V\right)$ .
\item As $cn\in M_{Limit}$ and as $cn_{i}\in\omega\left(M_{i}\right)$
($cn_{i}$ was a number) so must $cn\in_{Limit}\omega\left(M_{Limit}\right)$
be a number.
\end{enumerate}

\paragraph{Notation:}

As $cn\in_{Limit}\omega\left(M_{Limit}\right)$ is a number the notation
of $10\cdot cn$ is the number $cn$ times the number $10$.

\subsubsection{Formulas definition}

\subparagraph{Definition}

Let $\ensuremath{\phi\left(n,k\right)}$ be a two variable formula
(in $V$) which defines a function in $ZFC$ and holds the propriety
condition \ref{propriety-condition-definition} in $V$. Define that
following formula: 
\[
\zeta\left(\overset{\,o}{Z},n\right)=\exists x,k,c\in\omega\,\,\left(\begin{matrix}\phi\left(n,k\right)\,\wedge\,Run\left(x,n\right)=k\,\,\wedge\,x<\overset{\,o}{Z}\,\,\wedge\\
\left(\forall n',k'\in\omega\right)\left(\begin{matrix}\left(\phi\left(n',k'\right)\wedge Run\left(x,n'\right)=k'\right)\rightarrow\\
\left(Time\left(x,n'\right)=\left\lceil \log\left(n'+1\right)+\log\left(k'+1\right)+c\right\rceil \right)
\end{matrix}\right)
\end{matrix}\right)
\]

\subparagraph{Explanation:}

The formula $\zeta$ states that for given $\overset{\,o}{Z},n$ the
number $\overset{\,o}{Z}$ is large enough so that the value of the
function $\phi$ on the input $n$ is computed by a machine of size
smaller than $\overset{\,o}{Z}$ and which runs in linear time. Our
aim will be to prove that $M_{Limit}\vDash\left(\forall k\in\omega\right)\,\zeta\left(10\cdot cn,k\right)$
where $cn\in M_{Limit}$ was defined in \ref{subsec:Consistency-number-in-M_limit}. 

\subsubsection{Knowledge is consistent in $M_{Limit}$ \label{subsec:Knowladge is consistent in M_=00007BLimit=00007D }}

\paragraph{Theorem:}

Let $k'\in\omega\left(V\right)$ and let $a\in Lace_{k'}$ Assume
that:
\begin{itemize}
\item $M_{Limit}\vDash a\in\omega$ (where $a$ is used here as also the
appropriate equivalence class).
\item $cn\restriction_{<k'}\in Lace_{k'}$ 
\item both $cn$ and $a$ are in $var_{k'}$ 
\item The formula $\zeta\left(10\cdot cn,a\right)$ or its negation is in
$Formulas_{k+1}$ 
\item $M_{k';}\vDash a_{k'}<\left(G_{0}\left(M_{k'}\right)\right)_{3}$
\end{itemize}
Then the formulas $\zeta\left(10\cdot cn,a\right)$ (and not its negation)
is in $Formulas_{k+1}$ 

\paragraph{Proof:}

Recall that
\[
cn_{k'}=\max_{t\in\omega\left(M_{k'}\right)}\left\{ M_{k'}\vDash t-con\left(ZFC\right)\right\} 
\]
 and denote $\varphi$ the axiom ``$\exists y,c\in\omega$ s.t $y\leq10\cdot cn_{k'}$
and $y$ codes a TM which computes for every $n<\left(G_{0}\left(M_{k'}\right)\right)_{3}$
the value of $k''\in\omega$ s.t $\phi\left(n,k''\right)$ in $\left\lceil \log\left(n+1\right)+\log\left(k''+1\right)+c\right\rceil $
steps''. Let $j\in\omega\left(V\right)$. Recall that by the theorem
of knowledge in $M_{i}$ \ref{subsec:Knowledge in M_=00007Bi=00007D}
\[
M_{k'}\vDash j-con\left(ZFC\cup Formulas_{k'+1}\cup\varphi\right)
\]
 and notice that $\zeta\left(cn_{k'},a_{k'}\right)$ follows from
$\varphi$ (as such $y$ in $\varphi$ will be the $x$ in $\zeta\left(cn_{k'},a\right)$
that give us the coding of the TM). Therefore, as 
\[
M_{k'}\vDash j-con\left(ZFC\cup Formulas_{k'+1}\cup\varphi\right)
\]
 we know that $M_{k'}$ has a sequence of $j$ models, each a set
of its previous, and each holds $\varphi$. As $\zeta\left(cn_{k'},a_{k'}\right)$
follows from $\varphi$ it must be the case that these models also
hold $\zeta\left(cn_{k'},a_{k'}\right)$ and as such 
\[
M_{k'}\vDash j-con\left(ZFC\cup Formulas_{k'+1}\cup\zeta\left(cn_{k'},a_{k'}\right)\right)
\]
 holds. If $\neg\zeta\left(cn_{k'},a_{k'}\right)$ was in $Formulas_{k'+1}$
it would be the case that a formula and its negation are consistent,
which is a contradiction to $M_{k'}\vDash con\left(ZFC\right)$ (which
we assumed). Thus, $\zeta\left(cn_{k'},a_{k'}\right)$ was in $Formulas_{k'+1}$.
$\blacksquare$

\subsubsection{$M_{Limit}$ knows $\phi$ on every number \label{subsec:M_=00007BLimit=00007D knows =00005Cphi on every number}}

\paragraph{Lemma:}

Let $k'\in\omega\left(V\right)$ and let $a\in Lace_{k'}$ s.t for
all $i\in\omega\left(V\right)$ where $i>k'$, it holds $M_{k'}\vDash a_{k'}\in\omega\left(M_{k'}\right)$.

Then, for infinitely many $i\in\omega\left(V\right)$ where $i>k'$,
it holds that $M_{i}\vDash a_{i}<\left(G_{0}\left(M_{i}\right)\right)_{3}$

\paragraph{Proof:}

Assume that $k>k'$, we will show the existence of an $i>k$ such
that $M_{i}\vDash a_{i}<\left(G_{0}\left(M_{i}\right)\right)_{3}$.
Denote $E\subset P_{0}$ to be the set of partial function $f$ in
$V$ which agree with $G_{0}$ on $\left(M_{i}\right)_{i=1}^{k}$
and such that if we continue the construction according to $f$ we
will get $i>k$ such that $M_{i}\vDash a_{i}<\left(f\left(M_{i}\right)\right)_{3}$. 

Please note that, as the construction of $M_{i}$ depended on $f_{1},f_{2}$
alone and didn't depend on $f_{3}$ (see \ref{independence of the lattice from (G_=00007B0=00007D)_=00007B3=00007D})
any partial function $f'$ could be extended into a partial function
in $E$. As one could simply change the value of $\left(f\left(M_{i}\right)\right)_{3}$,
for $i$ not yet defined on, to be larger than $a_{i}$ (and as $M_{i}$
or $a_{i}$ didn't depend on this value, they remain the same). Thus
as every partial function in $P_{0}$ can be extended to a function
in $E$ and thus $E$ is a set in $V$ that is dense over $G_{0}\restriction_{\left(M_{i}\right)_{1\leq i\leq k}}$
and thus by \ref{subsec:Forcing} theorem 3 we get that $E\cap G_{0}\neq\emptyset$.
As such an $i>k$ exists such that $M_{i}\vDash a_{i}<\left(G_{0}\left(M_{i}\right)\right)_{3}$.
$\blacksquare$

\paragraph{Theorem:}

Let $a\in M_{Limit}$ s.t $M_{Limit}\vDash a\in\omega$ then $M_{Limit}\vDash\zeta\left(10\cdot cn,a\right)$ 

\paragraph{Proof:}

Recall that by theorem \ref{Stabilization in Formulas_=00007Bk=00007D}
exist a $k\in\omega\left(V\right)$ s.t either for all $k'\in\omega\left(V\right)$
s.t $k'>k$ the formula $\zeta\left(10\cdot cn_{k},a_{k}\right)$
is in $Formulas_{k'}$. Or, for all $k'\in\omega\left(V\right)$ s.t
$k'>k$ the formula $\neg\zeta\left(10\cdot cn_{k},a_{k}\right)$
is in $Formulas_{k'}$. 

By \ref{Appearance-of-formulas in Formulas_=00007Bk=00007D } the
formula $\zeta\left(10\cdot cn_{k},a_{k}\right)$ or its negation
must appear in $Formulas_{k}$ from a certain point on.

By the lemma above we know that for infinitely many $i\in\omega\left(V\right)$
where $i>k'$, it holds that $M_{i}\vDash a_{i}<\left(G_{0}\left(M_{i}\right)\right)_{3}$
(where $k'$ is large enough to include $cn$ and $a$ variables).

By \ref{subsec:Knowladge is consistent in M_=00007BLimit=00007D }
we know that when $M_{i}\vDash a_{i}<\left(G_{0}\left(M_{i}\right)\right)_{3}$
it must be the case that $\zeta\left(10\cdot cn,a\right)$ (and not
its negation) is in $Formulas_{k+1}$. As the formula $\zeta\left(10\cdot cn,a\right)$
appears infinitely many times in $Formulas_{i}$ for different $i$s
and as such a formula must stabilize we get that $\zeta\left(10\cdot cn_{k},a_{k}\right)$
appears on $Formulas_{k}$ from certain point on and as such by \ref{4.2.5 Truth value in M_=00007BLimit=00007D }
we get that $M_{Limit}\vDash\zeta\left(10\cdot cn,a\right)$. $\blacksquare$

\subsubsection{$M_{Limit}$ knows $\phi$ }

\paragraph{Theorem:}

It holds that $M_{Limit}\vDash\left(\forall k\in\omega\right)\,\zeta\left(10\cdot cn,k\right)$ 

\paragraph{Proof:}

Let $k\in M_{Limit}$ s.t $M_{Limit}\vDash\forall k\in\omega$ by
previous theorem in \ref{subsec:M_=00007BLimit=00007D knows =00005Cphi on every number}
we know that $M_{Limit}\vDash\zeta\left(10\cdot cn,k\right)$ and
as such $M_{Limit}\vDash\left(\forall k\in\omega\right)\,\zeta\left(10\cdot cn,k\right)$.
$\blacksquare$

\paragraph{Theorem:}

It holds that $M_{Limit}\vDash\left(\exists\overset{\,o}{Z}\in\omega\right)\left(\forall k\in\omega\right)\,\zeta\left(\overset{\,o}{Z},k\right)$ 

\paragraph{Proof:}

By the previous theorem the assertion hold with $\overset{\,o}{Z}=10\cdot cn$
and where $cn\in M_{Limit}$.$\blacksquare$ 

\section{Percolation theorems\label{sec:Percolation-theorems}}

\paragraph{Notation:}

Let $W$ be a $ZFC$ model with a worldly cardinal and let $\phi\left(x,y\right)$
be a a two variable formula (in $W$) which defines a function in
$ZFC$ and holds the propriety condition \ref{propriety-condition-definition}.
Denote the formula 
\[
know_{\phi}=\begin{matrix}\exists\overset{\,o}{Z}\in\omega\,\,\forall n\in\omega\,\,\exists x,k,c\in\omega\\
\left(\begin{matrix}\phi\left(n,k\right)\,\wedge\,Run\left(x,n\right)=k\,\,\wedge\,x<\overset{\,o}{Z}\,\,\wedge\\
\left(\forall n',k'\in\omega\right)\left(\begin{matrix}\left(\phi\left(n',k'\right)\wedge Run\left(x,n'\right)=k'\right)\rightarrow\\
\left(Time\left(x,n'\right)=\left\lceil \log\left(n'+1\right)+\log\left(k'+1\right)+c\right\rceil \right)
\end{matrix}\right)
\end{matrix}\right)
\end{matrix}
\]
 where $Run\left(x,n\right)$ denote that function that return the
value return by the TM numbered $x$ on input $n$ and $Time\left(x,n\right)$
return the number of steps done in the calculation of $x$ on input
$n$. 

\subsection{First percolation theorem \label{subsec:First-percolation-theorem}}

\paragraph{Theorem:}

Let $W$ be a $ZFC$ model with a worldly cardinal and let $\phi\left(x,y\right)$
be a a two variable formula (in $W$) which defines a function in
$ZFC$ and holds the propriety condition \ref{propriety-condition-definition}
. Then $know_{\phi}$ is consistent with $ZFC$.

\paragraph{Proof:}

Use \ref{subsec:Base-model-definitions} to construct $V,M_{1}$ then
use \ref{subsec:Forcing-POS-over-V} to define $P_{0}$ and then construct
the lattice \ref{subsec:Lattice-construction}. Define $M_{Limit}$
as in \ref{subsec:Limit-model-definition} and thus $M_{Limit}$ is
a $ZFC$ model which holds $know_{\phi}$. Thus a there is a model
of $ZFC$ that holds $know_{\phi}$, the formula $know_{\phi}$ is
consistent with $ZFC$.

\subsection{Second percolation theorem \label{subsec:Second-percolation-theorem}}

\paragraph{Theorem:}

Let $W$ be a $ZFC$ model with a worldly cardinal and let $\phi\left(x,y\right)$
be a a two variable formula (in $W$) which defines a function in
$ZFC$ and holds the propriety condition \ref{propriety-condition-definition}
. Then 
\[
W\vDash know_{\phi}
\]

\paragraph{Proof:}

Assume by negation that 
\[
W\not\vDash know_{\phi}
\]
. Then we can define a new axiom set 
\[
ZFC^{new}=ZFC\cup\neg know_{\phi}
\]

as $W$ is a model of $ZFC^{new}$ with a worldly cardinal we can
construct the same construction of \ref{subsec:First-percolation-theorem}
when using $ZFC^{new}$ instead of the regular $ZFC$. Thus we will
get a model of $ZFC^{new}$ , $M_{Limit}$ in which $know_{\phi}$
is held. But, $M_{Limit}$ also hold $\neg know_{\phi}$ as a model
of $ZFC^{new}$ which is a contradiction. Thus we get 
\[
W\vDash know_{\phi}
\]
 $\blacksquare$

\paragraph{Additional explanation:}

In light of the above proof. We'd like to offer additional explanation
on key points of the construction done under $ZFC^{new}$ of $M_{Limit}$: 
\begin{enumerate}
\item Recall that we start with $W$ a model with a worldly cardinal in
both cases in \ref{subsec:Base-model-definitions} and in \ref{subsec:Second-percolation-theorem}.
\item As the formula $\neg know_{\phi}$ is $\omega$ attributed it must
be held by any model with the same $\omega$. As such both $W$ and
$V_{k}^{W}$ hold $\omega-con\left(ZFC^{new}\right)$ by the theorem
in \ref{subsec:Base-model-definitions}. As such $V$ will hold it
as an elementary countable submodel.
\item $M_{1}$ is chosen the same as in \ref{subsec:M_1-construction} with
$ZFC^{new}$ replacing $ZFC$.
\item The forcing POS is done the same way and by the forcing theorem $V\left[G_{0}\right]$
is a $ZFC$ model. But as $\omega\left(V\right)=\omega\left(V\left[G_{0}\right]\right)$
the model $V\left[G_{0}\right]$ also hold $\omega-con\left(ZFC^{new}\right)$.
\item The rest of the construction of the lattice is done the same with
$ZFC^{new}$ replacing $ZFC$. And the construction of $M_{Limit}$
is done the same.
\item As each $M_{i}$ held $ZFC^{new}$ , by the same proof as in \ref{subsec:ZFC holds in M_=00007BLimit=00007D }
the model $M_{Limit}$ holds $ZFC^{new}$.
\item The same proofs of the lemma and the theorem of knowledge in $M_{i}$
of \ref{subsec:Knowledge in M_=00007Bi=00007D} applies with $ZFC^{new}$
replacing $ZFC$. 
\item The definition of the consistency number in \ref{subsec:Consistency-number-in-M_limit}
is done the same with $ZFC^{new}$ replacing $ZFC$. 
\item The same proofs of \ref{subsec:Knowladge is consistent in M_=00007BLimit=00007D }
and \ref{subsec:M_=00007BLimit=00007D knows =00005Cphi on every number}
holds the same with with $ZFC^{new}$ replacing $ZFC$. 
\end{enumerate}

\section{$P$ vs $NP\cap co-NP$ and other uses cases}

In this section we'll demonstrate the usage of section \ref{sec:Percolation-theorems}
in order to show that for a model $w$ with A worldly cardinal for
any language $L\in NP\cap co-NP$ provably, it is the case that $L\in P$.
In other words
\[
P=NP\cap co-NP
\]
assuming knowledge. Additionally, I'll give another use case of number
factoring. 

\subsection{Basic definition \label{subsec:Basic-definition NP-cap-co-NP} }

\paragraph{Definition:}

Let $\psi\left(x\right)$ be a formula with a single free variable.
The set 
\[
L_{\psi}=\left\{ x\in\omega\,|\,\psi\left(x\right)\right\} 
\]
 is called the language defined by $\psi$.

\paragraph{Definition: }

Let $V,\in_{V}$ be a $ZFC$ model. Let $\psi\left(x\right)$ be a
formula in $V$ with a single free variable. $\psi$ is said to define
a language in $NP\cap co-NP$ uniformly in $V$ if 

$m_{1},m_{2},k,s,t\in_{V}\omega\left(V\right)$ exist s.t the following
statements are \textbf{provable }from\textbf{ $ZFC$ }within\textbf{
$V$ }:
\begin{itemize}
\item $m_{1}$ and $m_{2}$ codes TM 
\begin{itemize}
\item $m_{1}$ will be the positive verifier and $m_{2}$ will be the negative
verifier.
\end{itemize}
\item the running time of $m_{1}$ on input $\left(x,y\right)$ and $m_{2}$
on input $\left(x,y\right)$ is at most $s\cdot\log{}^{k}x+t$ (for
any input $x,y$)
\item for any $x\in\omega$:
\begin{itemize}
\item if $\psi\left(x\right)$ then:
\begin{itemize}
\item exist $y\in\omega$ s.t $y\leq s\cdot x^{k}+t$ and $\left[m_{1}\right]\left(x,y\right)$
accepts.
\item for any $z\in\omega$ s.t $z\leq s\cdot x^{k}+t$ , $\left[m_{2}\right]\left(x,z\right)$
rejects.
\end{itemize}
\item if $\neg\psi\left(x\right)$ then:
\begin{itemize}
\item exist $y\in\omega$ s.t $y\leq s\cdot x^{k}+t$ and $\left[m_{2}\right]\left(x,y\right)$
accepts.
\item for any $z\in\omega$ s.t $z\leq s\cdot x^{k}+t$ , $\left[m_{1}\right]\left(x,z\right)$
rejects.
\end{itemize}
\end{itemize}
\end{itemize}

\paragraph{Remarks:}
\begin{enumerate}
\item The reader may want to think of the case $\omega\left(V\right)=\mathbb{N}$
at first.
\item the notation $\left[m\right]\left(x,y\right)$ is the run of $m$
on the inputs $x,y$
\item Please notice that we ask that the statements are \textbf{provable}.
meaning that all set models of $V$ will hold these conditions with
the ``same'' $m_{1},m_{2},k,s,t$ constants.
\item For the reader who is familiar with the $NP$ and $co-NP$ ``regular''
definitions, the above is the equivalent definitions using verifiers. 
\begin{enumerate}
\item The ``regular'' definition allows for different polynomials that
bounds the running time of $\left[m_{1}\right]$, $\left[m_{2}\right]$
and that bound the maximal length of the proofs. As we don't want
to use many constants for this (it'll make the notation even more
cumbersome). One can take the maximum of the powers and the maximum
of the free constant in order to get one bounding polynomial.
\end{enumerate}
\end{enumerate}

\subsection{construction's definitions and lemmas}

\paragraph{Definition:}

Let $V,\in_{V}$ be a $ZFC$ model. Let $\psi\left(x\right)$ be a
formula in $V$ that define a language in $NP\cap co-NP$ uniformly.
Let $m_{1},m_{2},s,k,t\in_{V}\omega\left(V\right)$ be the constants
as in \ref{subsec:Basic-definition NP-cap-co-NP}. Define the following
function 
\[
g\left(x\right)=\begin{cases}
1,y' & \psi\left(x\right)\\
0,y'' & \neg\psi\left(x\right)
\end{cases}
\]
 where $y'$ is the minimal number s.t $\left[m_{1}\right]\left(x,y'\right)$
accepts (positive proof) and $y''$ is the minimal number s.t $\left[m_{2}\right]\left(x,y''\right)$
accepts (negative proof). The numbers $0,1$ are just to indicate
that we code a negative \textbackslash{} positive poof.

Define $\phi_{\psi}\left(x,y\right)$ to be the formula $y=g\left(x\right)$
for the above $g\left(x\right)$.

\paragraph{Lemma:}

The formula $\phi_{\psi}\left(x,y\right)$ defines a function in $ZF$
and holds the propriety condition 

\paragraph{Proof:}

The fact that $\phi_{\psi}\left(x,y\right)$ defines a function is
obvious, as $\phi_{\psi}\left(x,y\right)$ it true iff $y=g\left(x\right)$
for a function $g$. One needs only prove that $\phi_{\psi}$ holds
the propriety condition. Let $V_{1}\stackrel[I]{f}{\rightarrow}V_{2}$
s.t $V_{2}$ is a set of $V_{1}$ and assume $x,y\in\omega\left(V_{1}\right)$
s.t $V_{1}\vDash\phi_{\psi}\left(x,y\right)$ in which case $y$ is
either $1,y'$ s.t $\left[m_{1}\right]\left(x,y'\right)$ accepts
or $0,y''$ s.t $\left[m_{2}\right]\left(x,y''\right)$ accepts. By
\ref{subsec:absoluteness-of-Turing} we know that if $\left[m_{1}\right]\left(x,y'\right)$
accepts then so $\left[f\left(m_{1}\right)\right]\left(f\left(x\right),f\left(y'\right)\right)$
and the same for $\left[m_{2}\right]\left(x,y''\right)$. Therefore,
if $x,y\in\omega\left(V_{1}\right)$ s.t $V_{1}\vDash\phi_{\psi}\left(x,y\right)$
then $V_{2}\vDash\phi_{\psi}\left(f\left(x\right),f\left(y\right)\right)$.
On the other hand if $x,y\in\omega\left(V_{1}\right)$ are s.t $V_{1}\vDash\neg\phi_{\psi}\left(x,y\right)$
then either 
\begin{enumerate}
\item $y$ isn't of the form $1,y'$ or $0,y''$ 
\item $y$ is of the form $1,y'$ but $\left[m_{1}\right]\left(x,y'\right)$
rejects
\item $y$ is of the form $0,y''$ but $\left[m_{2}\right]\left(x,y''\right)$
rejects
\end{enumerate}
In all cases $f\left(y\right)$ will hold:
\begin{enumerate}
\item $f\left(y\right)$ isn't of the form $1,f\left(y'\right)$ or $0,f\left(y''\right)$
\item $f\left(y\right)$ is of the form $1,f\left(y'\right)$ and $\left[f\left(m_{1}\right)\right]\left(f\left(x\right),f\left(y'\right)\right)$
rejects
\item $f\left(y\right)$ is of the form $0,f\left(y''\right)$ and $\left[f\left(m_{1}\right)\right]\left(f\left(x\right),f\left(y'\right)\right)$
rejects.
\end{enumerate}
and in all cases $V_{2}\vDash\neg\phi_{\psi}\left(f\left(x\right),f\left(y\right)\right)$.
$\blacksquare$ 

\subsection{Given knowledge, $P=NP\cap co-NP$ non-uniformly \label{Given knowledge, P=00003DNP=00005Ccap co-NP non-uniformly}}

\paragraph{Theorem:}

$M,\in_{M}$ be a model with a worldly cardinal. Let $\psi\left(x\right)$
be a formula in $M$ that define a language in $NP\cap co-NP$ uniformly.
Then a the language 
\[
\left\{ x\in_{M}\omega\left(M\right)\,|\,M\vDash\psi\left(x\right)\right\} \subset\omega\left(M\right)
\]
 is decidable in poly-logarithmic time in $M$

\paragraph{Proof:}

First we may assume that $k\left(n\right)$ s.t $\phi_{\psi}\left(n,k\left(n\right)\right)$.
Let $m_{1},m_{2},s,k,t\in_{M}\omega\left(M\right)$ be the constants
of $\psi\left(x\right)$ in $M$.

In $M$ it holds that $m_{1},m_{2}$ are positive and negative verifies
respectively and each run at at most $s\cdot\log^{k}x+t$ steps at
most, as this was proved from $ZFC$.

In $M$, by \ref{subsec:Second-percolation-theorem}, it holds 
\[
M\vDash\begin{matrix}\exists\overset{\,o}{Z}\in\omega\,\,\forall n\in\omega\,\,\exists x,k,c\in\omega\\
\left(\begin{matrix}\phi_{\psi}\left(n,k\right)\,\wedge\,Run\left(x,n\right)=k\,\,\wedge\,x<\overset{\,o}{Z}\,\,\wedge\\
\left(\forall n',k'\in\omega\right)\left(\begin{matrix}\left(\phi_{\psi}\left(n',k'\right)\wedge Run\left(x,n'\right)=k'\right)\rightarrow\\
\left(Time\left(x,n'\right)=\left\lceil \log\left(n'+1\right)+\log\left(k'+1\right)+c\right\rceil \right)
\end{matrix}\right)
\end{matrix}\right)
\end{matrix}
\]
 where $Run\left(x,n\right)$ denote that function that return the
value return by the TM numbered $x$ on input $n$ and $Time\left(x,n\right)$
return the number of steps done in the calculation of $x$ on input
$n$.

Define the following TM $T$ in $M$, given $n\in_{M}\omega$ :
\begin{itemize}
\item Until a $k$ s.t $\phi_{\psi}\left(n,k\right)$ is found:
\begin{enumerate}
\item run all TMs coded by numbers $y<_{M}\overset{\,o}{Z}$ one more step.
\item for each $y$, a TM that halted on the last step, check:
\begin{enumerate}
\item if the result of the calculation is of the form $1,y'$ for a number
$y'\in_{M}\omega$ which hold $y'\leq s\cdot n^{k}+t$ check if $\left[m_{1}\right]\left(n,y'\right)$.
If so. then $k=\left\langle 1,y'\right\rangle $ halt $T$ and return
true (as we've found that $\psi\left(n\right)$).
\item if the result of the calculation is of the form $0,y''$ for a number
$y''\in_{M}\omega$ which hold $y''\leq s\cdot n^{k}+t$ check if
$\left[m_{2}\right]\left(n,y''\right)$. If so. then $k=\left\langle 0,y''\right\rangle $
halt $T$ and return false (as we've found that $\neg\psi\left(n\right)$).
\end{enumerate}
\begin{itemize}
\item if both conditions (a) + (b) failed return to (1).
\end{itemize}
\end{enumerate}
\end{itemize}
The TM $T$ calculates $k$ s.t $\phi_{\psi}\left(n,k\right)$ (due
to the fact that $know_{\phi}$ holds in $M$). The question now is
its running time. 

\subparagraph{Running time analysis:}

The following running time analysis is done within $M$:
\begin{itemize}
\item each step of (1) takes $\overset{\,o}{Z}$ steps (assuming simulating
a TM one step takes also one step, if it takes $z'$ steps then step
(1) $z'\cdot\overset{\,o}{Z}$ steps). 
\item 2a takes at most $s\cdot\log^{k}y'+t$ steps. 
\item 2b takes at most $s\cdot\log^{k}y''+t$ steps. 
\item both $y',y''$ are bounded by $n^{k}+t$.
\item exist $y<_{M}\overset{\,o}{Z}$ that codes a TM which computes $k\left(n\right)$
in $\log\left(n+1\right)+\log\left(k+1\right)$$+c_{y}$ computing
steps. Therefore, the total number of iterations of step (1) is bounded
by 
\[
\overset{\,o}{Z}\cdot\left(\log\left(n+1\right)+\log\left(k+1\right)+\left(\max_{0\leq_{M}y\leq_{M}\overset{\,o}{Z}}c_{y}\right)\right)
\]
 when the maximum $\max_{0\leq_{M}y\leq_{M}\overset{\,o}{Z}}c_{y}$
is taken within $M$.
\item So, the total running time of this algorithm is poly-logarithmic time
bounded.
\end{itemize}

\subparagraph{Additional explanation of $\max_{0\protect\leq_{M}y\protect\leq_{M}\protect\overset{\,o}{Z}}c_{y}$:}

The reader may be baffled by the idea that $\max_{0\leq y\leq\overset{\,o}{Z}}c_{y}$
can be taken. However, $c_{y}$ can be defined within the model $M$,
as for each $y$ if $y$ codes a TM which computes $k\left(n\right)$
s.t $\phi\left(n,k\right)$ for every $n\leq t$ (for some $t$) in
$\log\left(n+1\right)+\log\left(k+1\right)$$+c_{y}$ steps we can
know the value $c_{y}$. If $y$ doesn't code such TM we can define
$c_{y}$ to be $0$. And thus the maximum $\max_{0\leq y\leq\overset{\,o}{Z}}c_{y}$
is a maximum of $\overset{\,o}{Z}$ numbers in $M$. 

$\blacksquare$

\paragraph{Remark:}

The term non-uniformly reminds us that the language we've started
with is provable in $NP\cap co-NP$. Thus, this assertion is different
from the assertion $P=NP\cap co-NP$ as in the latter, one must show
that for \uline{any} language in $NP\cap co-NP$ (regardless of
provability) that language is in $P$. Of course, we haven't showed
the latter statement in this paper.

\subsection{Number factorization}

In this subsection, I'll describe the usage of computational knowledge
to factor natural numbers. 

\subsubsection{Construction definition}

\paragraph{Definitions:}
\begin{itemize}
\item for $n,k\in\omega$ the notation ``$\mod\left(k,n\right)=0$'' denotes
that $k$ divides $n$. 
\item Let $n\in\omega$ 
\begin{itemize}
\item $n$ is called a composite if $n=a\cdot b$ for two numbers $a,b>1$.
\item let $n$ be a composite number. A number $b$ is called non-trivial
divisor of $n$ if $\mod\left(b,n\right)=0$ and $1<b<n$.
\end{itemize}
\end{itemize}

\paragraph{Definition:}

Let $V,\in_{V}$ be a $ZFC$ model. define the minimal divisor function
by

\[
g\left(n\right)=\min_{1<x\leq n}\left\{ x\,|\,\mod\left(x,n\right)=0\right\} 
\]
Define $\phi_{div}\left(x,y\right)$ to be the formula $y=g\left(x\right)$
for the above $g\left(x\right)$.

\subsubsection{Basic lemma}

\paragraph{Lemma:}

The formula $\phi_{\psi}\left(x,y\right)$ defines a function in $ZFC$
and holds the propriety condition 

\paragraph{Proof:}

The fact that $\phi_{div}\left(x,y\right)$ defines a function is
obvious, as $\phi_{div}\left(x,y\right)$ it true iff $y=g\left(x\right)$
for a function $g$. One needs only prove that $\phi_{div}$ holds
the propriety condition. Let $V_{1}\stackrel[I]{f}{\rightarrow}V_{2}$
s.t $V_{2}$ is a set of $V_{1}$ and assume $x,y\in\omega\left(V_{1}\right)$
s.t $V_{1}\vDash\phi_{div}\left(x,y\right)$ then $V_{1}\vDash\mod\left(y,x\right)=0$
and as such 
\[
V_{2}\vDash\mod\left(f\left(y\right),f\left(x\right)\right)=0
\]
 (as $f$ is arithmetic) and as $V_{1}\vDash\forall1<k<y\,\,\mod\left(k,x\right)\neq0$
it is the case that 
\[
V_{1}\vDash\forall1<k<f\left(y\right)\,\,\mod\left(f\left(k\right),f\left(x\right)\right)\neq0
\]
 (as $f$ is arithmetic). And thus $V_{2}\vDash\phi_{div}\left(f\left(x\right),f\left(y\right)\right)$
$\blacksquare$ 

\subsubsection{Given knowledge, factoring in in $P$}

\paragraph{Theorem:}

$M,\in_{M}$ be a model with a worldly cardinal then in $M$ then
exists, within $M$, $T_{1}$ a TM that given $n$ a composite number,
$T_{1}$ computes a non-trivial divisor.

\paragraph{Proof:}

In $M$, by \ref{subsec:Second-percolation-theorem}, it holds 
\[
M\vDash\begin{matrix}\exists\overset{\,o}{Z}\in\omega\,\,\forall n\in\omega\,\,\exists x,k,c\in\omega\\
\left(\begin{matrix}\phi_{div}\left(n,k\right)\,\wedge\,Run\left(x,n\right)=k\,\,\wedge\,x<\overset{\,o}{Z}\,\,\wedge\\
\left(\forall n',k'\in\omega\right)\left(\begin{matrix}\left(\phi_{div}\left(n',k'\right)\wedge Run\left(x,n'\right)=k'\right)\rightarrow\\
\left(Time\left(x,n'\right)=\left\lceil \log\left(n'+1\right)\right\rceil +\left\lceil \log\left(k'+1\right)\right\rceil +c\right)
\end{matrix}\right)
\end{matrix}\right)
\end{matrix}
\]
 where $Run\left(x,n\right)$ denote that function that return the
value return by the TM numbered $x$ on input $n$ and $Time\left(x,n\right)$
return the number of steps done in the calculation of $x$ on input
$n$.

Define the following TM $T_{1}$ in $M$, given $n\in_{M}\omega$
composite number to find a non trivial divisor of $n$:
\begin{itemize}
\item Until a $k$ s.t $1<k<n$ and $\mod\left(k,n\right)=0$ is found:
\begin{enumerate}
\item run all TMs coded by numbers $y<_{M}\overset{\,o}{Z}$ one more step.
\item for each $y$, a TM that halted on the last step, check if the result
of the calculation $z'$: 
\begin{enumerate}
\item check if $1<z'<n$ if so, go to (b) otherwise return to (1).
\item check if $\mod\left(z',n\right)=0$ if so, halt and return $z'$ otherwise
return to (1).
\end{enumerate}
\end{enumerate}
\end{itemize}
due to the fact that $know_{\phi_{div}}$ holds in $M$, we know that
for some $y<_{M}\overset{\,o}{Z}$ the answer of $g\left(n\right)$
will be given. We can't know that $g\left(n\right)$ will be the first
divisor to show up in the process but the TM $T_{1}$ calculates a
divisor, i.e a $k$ s.t $1<k<n$ and $\mod\left(k,n\right)=0$ . 

The question now is its running time. 

\subparagraph{Running time analysis:}

The following running time analysis is done within $M$:
\begin{itemize}
\item each step of (1) takes $\overset{\,o}{Z}$ steps (assuming simulating
a TM one step takes also one step, if it takes $z'$ steps then step
(1) $z'\cdot\overset{\,o}{Z}$ steps). 
\item step 2a is a simple comparison that takes at most $\log\left(n\right)$
\item step 2b takes at most poly-logarithmic time in the output of the machine.
as any candidate must be smaller than $n$ and division of two numbers
is done in poly-logarithmic time.
\item exist $y<_{M}\overset{\,o}{Z}$ that codes a TM which computes $g\left(n\right)$
in $\log\left(n+1\right)+\log\left(k+1\right)$$+c_{y}$ computing
steps. Therefore, the total number of iterations of step (1) is bounded
by 
\[
\overset{\,o}{Z}\cdot\left(\log\left(n+1\right)+\log\left(k+1\right)+\left(\max_{0\leq_{M}y\leq_{M}\overset{\,o}{Z}}c_{y}\right)\right)
\]
where $k$ is the divisor returned and when the maximum $\max_{0\leq_{M}y\leq_{M}\overset{\,o}{Z}}c_{y}$
is taken within $M$.
\item So, the total running time of this algorithm is poly-logarithmic time
bounded.
\end{itemize}

\paragraph{Theorem:}

$M,\in_{M}$ be a model with a worldly cardinal then in $M$ one can
factor numbers in poly-logarithmic time.

\paragraph{Proof:}

In $M$ use the following algorithm $T$, given $n\in\omega$:
\begin{enumerate}
\item verify that $n$ is a composite. If $n$ is prime, return $n$ as
farther factoring can't be done.
\item for $n$ a composite number run $T_{1}$ of the above theorem and
find a non-trivial factor $k$.
\item divide $n$ by $k$ and receive $\frac{n}{k}$, return the pair $n,\frac{n}{k}$.
\end{enumerate}
As by \cite{PRIMES is in P} we know that given a number n, testing
if n is prime or not is done in poly-logarithmic time. We know by
the theorem above that step 2 can be done in poly-logarithmic time
and we konw that step 3 can be done in poly-logarithmic time and as
such one can factor a number in poly-logarithmic time.

$\blacksquare$

\end{document}